\documentclass[letter,10pt,conference]{IEEEtran}

 \usepackage[bookmarks=false]{hyperref}
\ifCLASSINFOpdf
\else
\fi
\usepackage{relsize}
\usepackage{setspace}
\let\llncssubparagraph\subparagraph
\let\subparagraph\paragraph
\usepackage[compact]{titlesec}
\let\subparagraph\llncssubparagraph\usepackage{amsfonts}
\usepackage{mathrsfs}
\usepackage{bbm}
\usepackage{pifont}
\usepackage{amsfonts}
\usepackage{amsmath, amsthm}

\usepackage[top=0.625in,bottom=.75in,left=0.625in,right=0.625in]{geometry}
\usepackage{tabularx}
\usepackage{listings}
\usepackage{graphicx}
\usepackage{float}
\usepackage{amssymb}
\usepackage{array}
\usepackage{fancyhdr}
\usepackage{cases}
 \usepackage{supertabular}
\usepackage{url}
\usepackage{fancyhdr}

\usepackage{psfrag}
\usepackage{color}
\usepackage{bbding}
\newcommand{\bcap} {\hspace{2pt} \mathlarger{\cap}
\hspace{2pt}}

\newcommand{\bcup} {\hspace{2pt} \mathlarger{\cup}
\hspace{2pt}}

\newcommand {\C} {{\rm I\kern-5.5pt C}}

\newcommand{\bP}[1]{{\mathbb{P}}\left[{#1}\right]}

\def\centerhack#1{\hbox to 0pt{\hss\footnotesize #1\hss}}
\def\centerhackn#1{\hbox to 0pt{\hss #1\hss}}
\def\dchack#1{\vbox to 0pt{\vss{\hbox to 0pt{\hss#1\hss}}\vss}}

\AtBeginDocument{%
  \addtolength\abovedisplayskip{-0.3\baselineskip}%
  \addtolength\belowdisplayskip{-0.3\baselineskip}%
}

\newtheorem{lem}{Lemma}
\newtheorem{thm}{Theorem}

\newtheorem{cor}{Corollary}
\newtheorem{proposition}{Proposition}
\newtheorem*{proposition1.1}{Proposition 1.1}
\newtheorem*{proposition1.2}{Proposition 1.2}
\newtheorem*{proposition1.3}{Proposition 1.3}
\newtheorem*{proposition2.1}{Proposition 2.1}
\newtheorem*{proposition2.2}{Proposition 2.2}
\begin{document}

\title{Exact Analysis of $k$-Connectivity in\\Secure
Sensor Networks with Unreliable Links}

\author{ \IEEEauthorblockN{Jun Zhao}
\IEEEauthorblockA{CyLab and Dept.
of ECE \\
Carnegie Mellon University \\
Email: junzhao@cmu.edu} \and \IEEEauthorblockN{Osman Ya\u{g}an}
\IEEEauthorblockA{CyLab and Dept.
of ECE\\
Carnegie Mellon University \\
Email: oyagan@ece.cmu.edu} \and \IEEEauthorblockN{Virgil Gligor}
\IEEEauthorblockA{CyLab and Dept.
of ECE \\
Carnegie Mellon University \\
Email: gligor@cmu.edu}}

\maketitle \thispagestyle{plain} \pagestyle{plain}

%
%\author{\IEEEauthorblockN{Jun Zhao, Osman Ya\u{g}an and Virgil Gligor}
%\IEEEauthorblockA{\\CyLab and Dept.
%of ECE \\
%Carnegie Mellon University \\ \{junzhao, oyagan,
%virgil\}@andrew.cmu.edu}}

%\author{\IEEEauthorblockN{Jun Zhao, Osman Ya\u{g}an and Virgil Gligor}
%\IEEEauthorblockA{{\tiny~\vspace{-7pt}}\\CyLab and Dept.
%of ECE \\
%Carnegie Mellon University \\ \{junzhao, oyagan,
%virgil\}@andrew.cmu.edu\vspace{-7pt}}}

%\author{\IEEEauthorblockN{~}
%\IEEEauthorblockA{~\\
%~ \\ ~}}

\begin{abstract}
\boldmath The Eschenauer--Gligor (EG) random key predistribution
scheme has been widely recognized as a typical approach to secure
communications in wireless sensor networks (WSNs). However, there is
a lack of precise probability analysis on the reliable connectivity
of WSNs under the EG scheme. To address this,  we rigorously derive the asymptotically exact probability of $k$-connectivity
in
WSNs employing the EG scheme with unreliable links represented by independent on/off channels, where $k$-connectivity ensures that the
network remains connected despite the failure of any $(k-1)$ sensors
or links. Our analytical results are confirmed via numerical
experiments, and they provide precise guidelines for the design of
secure WSNs that exhibit a desired level of reliability against node
and link failures.

\end{abstract}

%\vspace{10pt}
%
%\begin{IEEEkeywords}
%key predistribution, node degree, random graph, random intersection
%graph, random key graph, security, topological properties, wireless
%sensor networks.
% \end{IEEEkeywords}


\begin{spacing}{.98}

\section{Introduction} \label{introduction}

The Eschenauer--Gligor (EG) random key predistribution scheme
\cite{virgil} has been widely regarded as a typical solution to
secure communications in wireless sensor networks (WSNs)
\cite{ISIT_RKGRGG,Krzywdzi,ryb3,zz,yagan_onoff,yagan,ISIT,ZhaoYaganGligor}.
The scheme operates as follows. In a WSN with $n$ sensors, before
deployment, each sensor is independently assigned $K_n$ distinct
keys which are selected \emph{uniformly at random} from a pool of
$P_n$ keys, where $K_n$ and $P_n$ are both functions of $n$. After
deployment, any two sensors can securely communicate over an
existing wireless link if and only if they share at least one key.

\iffalse

Random key predistribution schemes have been %extensively studied in
%the literature over the last decade and they are
 widely regarded as
appropriate solutions to secure communications in
%resource-constrained
 wireless sensor networks
\cite{r1,virgil,ryb3,zz,yagan_onoff,yagan,ZhaoYaganGligor,ISIT}. The
idea of randomly assigning cryptographic keys to sensors before
deployment has been introduced in the seminal work by Eschenauer and
Gligor \cite{virgil}. Their scheme, hereafter referred to as the EG
scheme, has received much interest
\cite{ryb3,zz,yagan_onoff,yagan,ZhaoYaganGligor,ISIT} and operates
as follows. In a WSN with $n$ sensors, prior to deployment, each
sensor is independently assigned $K_n$ distinct keys which are
selected \emph{uniformly at random} from a pool of $P_n$ keys, where
$K_n$ and $P_n$ are both functions of $n$. After deployment, any two
sensors can securely communicate over an existing wireless link if
and only if they share at least one key.
%There are several wireless communication models
%in the li
%Examples of physical link constraints include the reliability of the
%transmission channel and the distance between two sensors close
%enough for communication.

\fi

Wireless   links between nodes may become unavailable
 due to the presence of physical barriers
between nodes or because of harsh environmental conditions severely
impairing transmission. We model unreliable links as independent
channels, each being  {\em on} with probability $p_n $ or being {\em off} with probability $(1-p_n)$,
where $p_n$ is a function of $n$ for generality. Such on/off
channel model has been used in the context of secure WSNs
\cite{yagan_onoff,ZhaoYaganGligor,ISIT}, and is shown to well
approximate the disk model
\cite{ISIT_RKGRGG,Krzywdzi,yagan_onoff,ZhaoYaganGligor,ISIT}, where any two nodes need to be within a certain distance to
establish a wireless link in between.

Given the randomness involved in the EG key predistribution scheme,
and the unreliability of wireless links, there arises
a basic question as to how one can adjust the EG scheme parameters $K_n$
and $P_n$, and the link parameter $p_n$,
so that the resulting network is securely and
reliably connected. Reliability against the failure of sensors or links is
particularly important in WSN applications where sensors are
deployed in hostile environments (e.g., battlefield surveillance),
or, are unattended for long periods of time (e.g., environmental
monitoring), or, are used in life-critical applications (e.g.,
patient monitoring). To answer the question above, this paper presents the asymptotically exact
probability of $k$-connectivity in secure WSNs under the EG scheme with unreliable links. A network (or a
graph) is said to be $k$-connected if it remains connected despite
the deletion of any $(k-1)$ nodes or links. An equivalent definition is that each node can find at least $k$
internally node-disjoint paths to any other node. With $k = 1$, $k$-connectivity simply means connectivity.

%$k$-connectivity and the property that
%the minimum degree is at least $k$ are particularly beneficial in
%secure WSNs with autonomously operating sensors. A graph property
%weaker than and often related to $k$-connectivity is that the
%minimum (node) degree of the graph is at least $k$; i.e., each node
%is directly connected to no less than $k$ other nodes, where the
%minimum node degree refers to the minimum among the numbers of
%neighbors that nodes have.
%
%Under the on/off channel model, we perform a precise mathematical
%analysis on $k$-connectivity for arbitrary $k$ and the minimum node
%degree in WSNs that employ the EG scheme. Given the scheme
%parameters $K_n$ and $P_n$, and the probability $p_n$ of a wireless
%channel being \emph{on}, we derive the asymptotically exact
%probability of $k$-connectivity.

Our result on the asymptotically exact probability of
$k$-connectivity complements a zero-one law established in our prior work
\cite{ZhaoYaganGligor,ISIT}, and is significant to obtain a precise
understanding of the connectivity behavior of secure WSNs. First, with the
zero-one law, one is only provided with design choices which lead to networks
that are
$k$-connected with high probability or to that are not
$k$-connected with high probability, where an event happens ``with
high probability'' if its probability asymptotically converges to 1.
Given the trade-offs involved between connectivity, security and
memory load \cite{virgil,yagan_onoff}, it would be more useful to
have a complete picture by obtaining the asymptotically exact
probability of $k$-connectivity. In
addition, there may be situations where the network designer is
interested in having a guaranteed level of $k$-connectivity (one-laws
would provide conditions for that) but may also be interested in
having some level of $k$-connectivity without such guarantees
(one-laws would fall short in providing this). Our result fills this
gap. Finally, it is not possible to determine the width of the phase
transition from zero-one laws; the width of the phase transition is
often calculated by the difference in parameters that it takes to
increase the probability of $k$-connectivity from $\epsilon$ to
$(1-\epsilon)$, for some $\epsilon<0.5$. In other words, it is not
clear from zero-one laws how sensitive the probability of
$k$-connectivity is to the variations in the EG scheme parameters $K_n$
and $P_n$, and the link parameter $p_n$. By providing the
asymptotically exact probability of $k$-connectivity, our findings
provide a clear picture of these intricate relationships.

%Our approach is based on modeling the WSN by an intersection of two random
%graphs models; one induced by the EG scheme, and another by the on/off communication model.
%This renders the analysis challenging
%due to the intertwining of the two distinct types of random graphs
%\cite{yagan_onoff}.

%in secure WSNs
%employing the EG key predistribution scheme under the
%\emph{on}/\emph{off} channel model as the physical link constraint
%comprising independent channels which are either \emph{on} or
%\emph{off}.

The rest of the paper is organized as follows. We describe the system model in Section
\ref{sec:SystemModel}. Section \ref{sec:main:res}
presents the main results as Theorem 1, which is established in Section
\ref{sec_est}. In Section \ref{sec:expe}, we
present numerical experiments that confirm our analytical findings. Afterwards,
 Section \ref{related} surveys related work, and Section \ref{sec:Conclusion} concludes the paper. The Appendix presents a few
useful lemmas and their proofs.

 %
%
%
%offers
%
%discuss

\section{System Model} \vspace{-3pt}
\label{sec:SystemModel}

We now explain the system model. Consider a
WSN with $n$ sensors operating under the EG
scheme and with wireless links modeled by independent {on/off}
channels. Let a node set $\mathcal {V} = \{v_1, v_2, \ldots, v_n \}$ represent the $n$ sensors. According to the EG scheme,
each node $v_i \in \mathcal {V}$ is independently assigned a set
(denoted by $S_i$) of $K_n$ distinct cryptographic keys, which are
selected {\em uniformly at random} from a key pool of $P_n$ keys.
Any pair of nodes can then secure an existing communication link as
long as they have at least one key in common.

The EG scheme results in a \emph{random key
graph}
 \cite{r1,ryb3,yagan}, also known as a \emph{uniform random intersection
graph}. This graph denoted by
$G(n,K_n,P_n)$ is defined on the
node set $\mathcal{V}$ such that any two distinct nodes $v_i$ and
$v_j$ have an edge in between, an event denoted by $\Gamma_{ij}$, if
and only if they share at least one key. Thus, the event $\Gamma_{ij} $ means $ \big( S_{i} \bcap S_{j}\neq
\emptyset\big) $.

Under the {on/off} channel model for unreliable links, each wireless link is independently
being  {\em on} with probability $p_n $ or being {\em off} with probability $(1-p_n)$.
Defining ${C}_{i j}$ as the event that the channel between $v_i$ and
$v_j$ is {\em on}, we have $\bP{C_{ij}} = p_n$, with $\mathbb{P}[A]$
throughout the paper meaning the probability that event $A$ happens.
The {on/off} channel model induces an \emph{Erd\H{o}s-R\'enyi graph}
$G(n, p_n)$ \cite{citeulike:4012374} defined on the node set
$\mathcal{V}$ such that $v_i$ and $v_j$ have an edge in between if
$C_{ij}$ takes place.

 Finally, we denote by
$\mathbb{G}\iffalse_{on}\fi (n, K_n, P_n, p_n)$ the underlying graph
of the $n$-node WSN under the EG scheme with unreliable links. We often write $\mathbb{G}$ rather than
$\mathbb{G}(n, K_n, P_n, p_n)$ for brevity. Graph
$\mathbb{G}\iffalse_{on}\fi$ is defined on the node set
$\mathcal{V}$ such that there exists an edge between nodes $v_i$ and
$v_j$ if events $\Gamma_{ij}$ and $C_{ij}$ happen at the same time.
We set event $E_{ij} : = \Gamma_{ij} \cap C_{ij}$ and also write
$E_{ij} $ as $E_{v_i v_j} $ when necessary.
%\begin{equation}
%E_{ij} = \Gamma_{ij} \cap C_{ij},\textrm{ for }1\leq i < j \leq n,
%\label{eq:E_is_K_cap_C_oy}
%\end{equation}
 It is clear that $\mathbb{G}\iffalse_{on}\fi$ is the intersection of
$G(n, K_n, P_n)$ and $G(n, p_n)$; i.e.,
\begin{equation}
\mathbb{G}\iffalse_{on}\fi  = G(n, K_n, P_n) \cap G(n, p_n).
\vspace{-2pt} \label{graphinter}
%\label{eq:G_on_is_RKG_cap_ER_oy}
\end{equation}

We define $s_n $ as the probability that two distinct nodes share
at least one key and ${q_n} $ as the probability that two distinct
nodes have an edge in between in graph $\mathbb{G}$. Clearly, $s_n
$ and $q_n$ both depend on $K_n$ and $P_n$, while $q_n$ depends also
on $p_n$. As shown in previous work \cite{r1,ryb3,yagan}, $s_n$ is
determined through
\begin{align}
s_n & =  \mathbb{P} [\Gamma_{i j} ]  = \begin{cases} 1-
\binom{P_n- K_n}{K_n} \big/ \binom{P_n}{K_n} , &\textrm{if }P_n > 2
K_n , \\ 1, &\textrm{if }P_n \leq 2 K_n . \end{cases} \nonumber
\vspace{-2pt}
\end{align}
%and
%\begin{align}
%s_n & =
%\end{align}
Then by the independence of ${C}_{i j} $ and $ \Gamma_{i j} $, we
have
\begin{align}
{q_n}   & =  \mathbb{P} [E_{i j} ]  =  \mathbb{P} [{C}_{i j} ] \cdot
\mathbb{P} [\Gamma_{i j} ] =  p_n\iffalse_{on}\fi \cdot s_n
 \label{qnexprnewexpr1}
\\ & =
\begin{cases} p_n\iffalse_{on}\fi \cdot \big[1- \binom{P_n- K_n}{K_n}
 \big/ \binom{P_n}{K_n}\big], &\textrm{if }P_n > 2 K_n ,
\\ p_n, &\textrm{if }P_n \leq 2 K_n .
\end{cases} \label{qnexprnewexpr}
\end{align}

\section{The Main Results} \label{sec:main:res}

We present the main results below. Throughout the
paper, $k$ is a positive integer and does not scale with $n$, and $e$ is the base of
the natural logarithm function, $\ln$. We use the standard
asymptotic notation $o(\cdot), O(\cdot), \omega(\cdot),
\Omega(\cdot),\Theta(\cdot)$ and $ \sim$; in particular, for two
positive sequences $a_n$ and $b_n$, the relation $a_n \sim b_n$ means
$\lim_{n \to
  \infty} {a_n}/{b_n}=1$.%\footnote{Specifically, given two positive
%functions $f(n)$ and $g(n)$,
%\begin{enumerate}
%  \item $f(n) = o \left(g(n)\right)$ means $\lim_{n \to
%  \infty}\frac{f(n)}{g(n)}=0$.
%  \item $f(n) = O \left(g(n)\right)$ means that there exists a positive
%  constant $c_1$ such that $f(n) \leq c_1 g(n)$ for all $n$ sufficiently large.
%  \item $f(n) = \Omega \left(g(n)\right)$ means that there exists a positive
%  constant $c_2$ such that $f(n) \geq c_2 g(n)$ for all $n$ sufficiently
%  large.
%  \item $f(n) = \Theta \left(g(n)\right)$ means that there exist positive
%  constants $c_3$ and $c_4$ with $c_3\leq c_4$ such that
%  $c_3 g(n) \leq f(n) \leq c_4 g(n)$ for all $n$ sufficiently large.
%
%  \item $f(n) \sim g(n)$ signifies $\lim_{n \to
%  \infty}\frac{f(n)}{g(n)}=1$; namely, $f(n)$
%  and $g(n)$ are asymptotically equivalent.
%  \end{enumerate}
%}.
%
%Consider a positive integer $k$ and scalings $K \hspace{-3pt}:
%\mathbb{N}_0 \rightarrow \mathbb{N}_0,P \hspace{-3pt}: \mathbb{N}_0
%\rightarrow \mathbb{N}_0$ and $p \hspace{-2pt}: \mathbb{N}_0
%\rightarrow (0,1]$, with $P_n \geq 3K_n $ for all $n$ sufficiently
%large. Let the sequence $\alpha\hspace{-2pt}:
%\hspace{1pt}\mathbb{N}_0 \rightarrow \mathbb{R}$ be defined through

\begin{thm}\label{THM1}

For graph $\mathbb{G}(n, K_n, P_n, p_n)$ under $P_n = \Omega (n)$
and $\frac{K_n}{P_n} = o (1)$, with $q_n$ denoting the edge
probability and a sequence $\alpha_n$ defined through\vspace{-1.5pt}
\begin{align}
q_n & = {\frac{\ln  n + {(k-1)} \ln \ln n +
{\alpha_n}}{n}},\vspace{-1.5pt} \label{thm_eq_pe}
\end{align}
if $\lim_{n \to \infty} \alpha_n = \alpha ^* \in (-\infty, \infty)$,
then as $n \to \infty$,\vspace{-2pt}
\begin{align}
 \mathbb{P} \left[\hspace{2pt}\textrm{Graph }\mathbb{G}(n, K_n, P_n, p_n)\iffalse_{on}\fi \textrm{
is $k$-connected}.\hspace{2pt}\right] & \to
  e^{- \frac{e^{-\alpha ^*}}{(k-1)!}} .\vspace{-2pt} \nonumber %\label{eqn_kconn}
 \end{align}

\end{thm}

%We detail
%the proof of Theorem \ref{THM1} in the next section.
Theorem \ref{THM1} provides the asymptotically exact probability of $k$-connectivity in graph ${\mathbb{G}} $. Its proof is given in the next section. 
 From (\ref{qnexprnewexpr}), for all $n$ sufficiently large, under
$P_n > 2 K_n$ which is clearly implied by
 the condition $\frac{K_n}{P_n} = o (1)$, the edge probability $q_n$ in graph $\mathbb{G}$ is given by the\vspace{1pt} expression $p_n\iffalse_{on}\fi \cdot \big[1- \binom{P_n- K_n}{K_n}
 \big/ \binom{P_n}{K_n}\big]$. With a much simpler approximation $p_n \cdot \frac{{K_n}^2}{P_n} $ for $q_n$, we present below a corollary of Theorem \ref{THM1}.\vspace{-2pt}

 \begin{cor}\label{COR1}

For graph $\mathbb{G}(n, K_n, P_n, p_n)$ under $P_n = \Omega (n)$
and $\frac{{K_n}^2}{P_n} = o \big(\frac{1}{\ln n}\big)$, with a sequence $\beta_n$ defined through\vspace{-1.5pt}
\begin{align}
\textstyle{p_n \cdot \frac{{K_n}^2}{P_n}}   & = { \frac{\ln  n +
{(k-1)} \ln \ln n + {\beta_n}}{n}},\vspace{-1.5pt}
\label{thm_eq_pebeta}
\end{align}
if $\lim_{n \to \infty} \beta_n = \beta ^* \in (-\infty, \infty)$,
then as $n \to \infty$,\vspace{-2pt}
\begin{align}
 \mathbb{P} \left[\hspace{2pt}\textrm{Graph }\mathbb{G}(n, K_n, P_n, p_n)\iffalse_{on}\fi \textrm{
is $k$-connected}.\hspace{2pt}\right] & \to
  e^{- \frac{e^{-\beta ^*}}{(k-1)!}} . \vspace{-2pt}\nonumber %\label{eqn_kconn}
 \end{align}

\end{cor}

Setting
$p_n= 1$ in Theorem
\ref{THM1} and Corollary \ref{COR1}, we obtain the corresponding results for random key graph
$G(n,K_n,P_n)$ in view of (\ref{graphinter}).
Furthermore, we can use monotonicity arguments \cite{ZhaoYaganGligor} to derive
the zero-one laws for $k$-connectivity in graph $\mathbb{G}$. Specifically, under the conditions of Theorem
\ref{THM1} (resp., Corollary \ref{COR1}), graph $\mathbb{G}$ is $k$-connected with high probability if $\textstyle{\lim_{n \to \infty}}\alpha_n = \infty$ (resp., $\textstyle{\lim_{n \to \infty}}\beta_n = \infty$), and is not $k$-connected with high probability if $\textstyle{\lim_{n \to \infty}}\alpha_n = -  \infty$ (resp., $\textstyle{\lim_{n \to \infty}}\beta_n = - \infty$). The arguments are straightforward from our work \cite{ZhaoYaganGligor} and are omitted here due to space limitation.

%Namely, we have
%\begin{eqnarray}
%\textstyle{\lim_{n \to \infty}} \mathbb{P} \left[\mathbb{G}\iffalse_{on}\fi
%\textrm{ is $k$-connected}\hspace{2pt}\right]  =
%\begin{cases}
%1,  & \textrm{if} ~  \textstyle{\lim_{n \to \infty}}\alpha_n = \infty,   \\
% 0,   &   \textrm{if}~  \textstyle{\lim_{n \to \infty}}\alpha_n = - \infty.
%\end{cases}
%\nonumber
%\end{eqnarray}

Before establishing Corollary \ref{COR1} using Theorem \ref{THM1}, we explain the practicality of the conditions in Theorem
\ref{THM1} and Corollary \ref{COR1}: $P_n = \Omega (n)$, $\frac{K_n}{P_n} = o (1)$
and $\frac{{K_n}^2}{P_n} = o \big(\frac{1}{\ln n}\big)$.
 First, the condition $ P_n = \Omega(n)$
indicates that the key pool size $P_n$ should grow at least
linearly with $n$, which holds in practice
\cite{virgil,yagan,yagan_onoff}. Second, the condtions $\frac{K_n}{P_n} = o (1)$
and $\frac{{K_n}^2}{P_n} = o \big(\frac{1}{\ln n}\big)$ (note that the latter implies the former) are also practical in secure
sensor network applications since $P_n$ is expected to be several orders of magnitude larger than $K_n$ \cite{virgil,yagan,yagan_onoff}.

 We now prove Corollary \ref{COR1} using Theorem \ref{THM1}. We have the conditions of Corollary \ref{COR1}: $P_n = \Omega (n)$,
 $\frac{{K_n}^2}{P_n} = o \big(\frac{1}{\ln n}\big)$, and (\ref{thm_eq_pebeta}) with $\lim_{n \to \infty} \beta_n = \beta ^* \in (-\infty, \infty)$. First, it is clear that $\beta_n = \beta ^* \pm o(1) $. Under $\frac{{K_n}^2}{P_n} = o \big(\frac{1}{\ln n}\big) = o(1)$, from \cite[Lemma 8]{ZhaoYaganGligor}, it holds that $s_n = \frac{{K_n}^2}{P_n}  \cdot \big[1\pm O\big(\frac{{K_n}^2}{P_n} \big)\big] $. In view of the above, we obtain from (\ref{qnexprnewexpr1}) and (\ref{thm_eq_pebeta}) that\vspace{-1pt}
 \begin{align}
q_n     & = p_n \cdot s_n = \textstyle{p_n \cdot \frac{{K_n}^2}{P_n}  \cdot \big[1\pm O\big(\frac{{K_n}^2}{P_n} \big)\big]} \nonumber \\
& =\textstyle{ \frac{\ln  n + {(k-1)} \ln \ln n + {\beta_n}}{n} \cdot \big[1\pm o \big(\frac{1}{\ln n}\big)\big]}   \nonumber \\
& = \textstyle{\frac{\ln  n + {(k-1)} \ln \ln n +  \beta ^* \pm o(1)}{n}} . \vspace{-1pt} \label{qnbetao1}
\end{align}
With $\alpha_n$ defined by (\ref{thm_eq_pe}), we use (\ref{qnbetao1}) to derive
$\alpha_n = \beta ^* \pm o(1)$, which yields that $ \alpha ^*$ denoting $\lim_{n \to \infty} \alpha_n $ equals $\beta ^*$. Then in view of $ \alpha ^* = \beta ^*$ and that the conditions of Theorem \ref{THM1} all hold given the conditions of Corollary \ref{COR1} (note that  $\frac{{K_n}^2}{P_n} = o \big(\frac{1}{\ln n}\big)$ implies $\frac{K_n}{P_n} = o (1)$), Corollary \ref{COR1} follows from Theorem \ref{THM1}.

\section{Establishing Theorem \ref{THM1}} \label{sec_est}

%We will write $\mathbb{G}(n, K_n, P_n, p_n)$ as $\mathbb{G}$ at some
%places for notation brevity.

For any graph, $k$-connectivity implies
that its minimum degree is at least $k$, while the other way does
not hold since a graph may have isolated components, each of which
is $k$-connected within itself. However, for random graph
$\mathbb{G}(n, K_n, P_n, p_n)$, as given by Lemma \ref{lem-kcon-mnd}
below, we have shown it is unlikely under certain conditions that
$\mathbb{G}(n, K_n, P_n, p_n)$ is not $k$-connected but has a
minimum degree at least $k$.
\begin{lem}[\hspace{-.2pt}{\cite[Section IX]{ZhaoYaganGligor}}]
\label{lem-kcon-mnd}

For graph $\mathbb{G}(n, K_n, P_n, p_n)$ under $P_n = \Omega (n)$,
$\frac{K_n}{P_n} = o (1)$ and $q_n = o(1)$, it holds that
\begin{align}
 \mathbb{P}
 \left[\begin{array}{c}\textrm{Graph }\mathbb{G} \iffalse_{on}\fi \textrm{
is not $k$-connected},\\\textrm{but has a minimum degree at
least $k$}.\end{array}\right] & = o(1) . \nonumber %\label{eqn_kconn}
 \end{align}

\end{lem}

We show that the conditions in Lemma \ref{lem-kcon-mnd} all hold given
the conditions of Theorem \ref{THM1}: $P_n = \Omega (n)$, $\frac{K_n}{P_n} = o (1)$ and $q_n = \frac{\ln
n + {(k-1)} \ln \ln n + {\alpha_n}}{n}$ with $\lim_{n \to \infty}
\alpha_n = \alpha ^* \in (-\infty, \infty)$.
To see this, we only need to prove $q_n = o(1)$ needed in Lemma \ref{lem-kcon-mnd}  follows from the conditions of Theorem
\ref{THM1}. Clearly, it holds that $|\alpha_n|=O(1)$ from $\lim_{n \to
\infty} \alpha_n = \alpha ^* \in (-\infty, \infty)$. Then in view of
$|\alpha_n|=O(1)$ and the fact that $k$ does not scale with $n$, we
obtain from (\ref{thm_eq_pe}) that
\begin{align}
q_n & \sim  \frac{\ln n}{n},\label{eq_pe_lnnn}
\end{align}
which clearly implies $q_n = o(1)$.

From Lemma \ref{lem-kcon-mnd} and
\begin{align}
 & \mathbb{P} \left[\hspace{2pt}\textrm{Graph }\mathbb{G} \iffalse_{on}\fi \textrm{
is $k$-connected}.\hspace{2pt}\right] \nonumber \\ & = \mathbb{P}
\left[\hspace{2pt}\textrm{Graph }\mathbb{G} \iffalse_{on}\fi
\textrm{ has a minimum degree at least $k$}.\hspace{2pt}\right]
\nonumber \\ & \quad - \mathbb{P}
 \left[\begin{array}{c}\textrm{Graph }\mathbb{G} \iffalse_{on}\fi \textrm{
is not $k$-connected},\\\textrm{but has a minimum degree at least
$k$}.\end{array}\right], \nonumber
 \end{align}
 Theorem \ref{THM1} on $k$-connectivity of $\mathbb{G}$
will be proved once we demonstrate Lemma \ref{lem-mnd} below on
the minimum degree of $\mathbb{G}$.

%
%Theorem \ref{THM1} on $k$-connectivity will be proved if we show
%under conditions of Theorem \ref{THM1} that
%\begin{align}
% \mathbb{P} \left[\hspace{2pt}\textrm{Graph }\mathbb{G} \iffalse_{on}\fi \textrm{
%has a minimum degree at least $k$}.\hspace{2pt}\right] & \to
%  e^{- \frac{e^{-\alpha ^*}}{(k-1)!}} . \label{mndres}
% \end{align}

%
%We will establish (\ref{mndres}) under $K_n = \omega (1)$ and $P_n
%\geq 3K_n $ for all $n$ sufficiently large, which are more general
%conditions than the conditions of Theorem \ref{THM1} as explained
%below. By \cite[Lemma 7]{ZhaoYaganGligor}, $P_n = \Omega (n)$ and .
%Then the proof of Theorem \ref{THM1} will be completed once we
%demonstrate Lemma \ref{lem-mnd} below.

\begin{lem} \label{lem-mnd} %\hspace{2pt}

Under the conditions of Theorem \ref{THM1},
it holds that $\lim_{n\to\infty} \mathbb{P} [\mathbb{G} \iffalse_{on}\fi \textrm{
has a minimum degree at least $k$}.]=
  e^{- \frac{e^{-\alpha ^*}}{(k-1)!}} .$
%\begin{align}
% \mathbb{P} \left[\hspace{2pt}\textrm{Graph }\mathbb{G} \iffalse_{on}\fi \textrm{
%has a minimum degree at least $k$}.\hspace{2pt}\right] & \to
%  e^{- \frac{e^{-\alpha ^*}}{(k-1)!}} . \nonumber
% \end{align}
\end{lem}

To prove Lemma \ref{lem-mnd}, we first show that the number of nodes
in $\mathbb{G}\iffalse_{on}\fi$ with a certain degree converges in
distribution to a Poisson random variable. With $\phi_h$ denoting
the number of nodes with degree $h$ in $\mathbb{G}\iffalse_{on}\fi$,
$h = 0,1, \ldots$, we use the method of moments to prove that
$\phi_h $ asymptotically follows a Poisson distribution with mean
$\lambda_h$. Specifically, from %\cite[Theorem 2.13]{2008asymptotic}
 \cite[Theorem 7]{ChihWeiYi_Globecom07}, it follows for any
integers $h \geq 0$ and $\ell \geq 0$ that
\begin{align}
 \mathbb{P}[\phi_h = \ell]
 & \sim (\ell !)^{-1}{\lambda_h} ^{\ell}e^{-\lambda_h},
 \label{eqn_phihell}
 \end{align}
 since $\mathbb{P} [\textrm{Nodes }v_{1}, v_{2}, \ldots, v_{m}\textrm{ all have
degree }h]    \sim {\lambda_h}^m / n^m$, which is shown by Lemma \ref{LEM1} below with
\begin{align}
 \lambda_h & = n (h!)^{-1}(n q_n)^h e^{-n
q_n}. \label{eqn_labmdah}
 \end{align}

%\begin{lem} \label{lem_pos_sum}
%
%For any constant integers $m \geq 1$ and $h \geq 0$, it holds that
%\begin{align}
% &  \sum_{\begin{subarray}{c} (i_1, i_2, \ldots, i_m):
%   \\
%1 \leq i_1 < i_2 < \ldots < i_m \leq n
%\end{subarray}} \mathbb{P} [\textrm{Nodes }v_{i_1}, v_{i_2}, \ldots, v_{i_m}\textrm{ have degree
%}h]
%\nonumber  \\
% & \hspace{215pt} \sim \frac{{\lambda_h}^m}{m!} , \nonumber
%\end{align}
%where $\lambda_h$ is given by (\ref{eqn_labmdah}).
%
%\end{lem}

\begin{lem} \label{LEM1}

For graph $\mathbb{G}\iffalse_{on}\fi$ under the conditions of Theorem \ref{THM1},
 $\mathbb{P}   [v_{1}, \hspace{-1pt}v_{2},\hspace{-1pt} \ldots, \hspace{-1pt}v_{m}\hspace{3pt}\textrm{all have
degree }h]  \hspace{-1pt} \sim \hspace{-1pt}  (h!)^{-m}\hspace{-1pt}  (n q_n)^{hm}\hspace{-1pt} e^{-m n q_n}$ holds for any integers $m \geq 1$ and $h
\geq 0$. % we have
%\begin{align}
% &  \mathbb{P} [\textrm{Nodes }v_{1}, v_{2}, \ldots, v_{m}\textrm{ have
%degree }h] \nonumber  \\
% & \quad\quad\quad\quad\quad\quad\quad\quad\quad\quad~ \sim  (h!)^{-m}  (n q_n)^{hm} e^{-m n q_n}.\nonumber
%\end{align}

\end{lem}

As explained above, Lemma \ref{LEM1} shows (\ref{eqn_phihell}) with
$\lambda_h$ given by (\ref{eqn_labmdah}). Then the proof of Lemma
\ref{lem-mnd} will be completed once we establish Lemma \ref{LEM1}
and the result that (\ref{eqn_phihell}) implies Lemma \ref{lem-mnd}.
Below we will demonstrate that (\ref{eqn_phihell}) implies Lemma
\ref{lem-mnd}, and then detail the proof of Lemma \ref{LEM1}.

%We first introduce useful results on $q_n$, the edge probability of
%graph $\mathbb{G}\iffalse_{on}\fi$.

\subsection{Proving that (\ref{eqn_phihell}) implies Lemma
\ref{lem-mnd}}

Recall that $\phi_h$ denotes the number of nodes with degree $h$ in
graph $\mathbb{G}\iffalse_{on}\fi$. With $\delta$ defined as the
minimum
 degree of graph $\mathbb{G}\iffalse_{on}\fi$, then the event
$(\delta \geq k)$ is the same as $ \bigcap_{h=0}^{k-1}
(\phi_h = 0) $ (i.e., the event that no node has a degree falling in
$\{0,1,\ldots, k-1\}$). Hence, we obtain
\begin{align}
 \mathbb{P}[\delta \geq k] & = \mathbb{P}\bigg[\bigcap_{h=0}^{k-1}
(\phi_h = 0)\bigg] \leq   \mathbb{P}[\phi_{k-1} = 0];
\label{eqn_1mindel2}
 \end{align}
and by the union bound, it holds that
\begin{align}
\mathbb{P}[\delta \geq k] & = \mathbb{P}\bigg[ ~(\phi_{k-1} = 0)
\hspace{2pt} \mathlarger{\cap} \hspace{2pt}
\overline{\bigg(\bigcup_{h=0}^{k-2} (\phi_h \neq 0) \bigg)}
 ~\bigg]  \nonumber  \\
 & \geq \mathbb{P}[\phi_{k-1} = 0] -
  \sum_{h=0}^{k-2}
\mathbb{P}[  \phi_{h} \neq 0] .
 \label{eqn_1min}
 \end{align}
To use (\ref{eqn_1mindel2}) and (\ref{eqn_1min}), we compute
$\mathbb{P}[ \phi_{h} \neq 0]$ given (\ref{eqn_phihell}) and thus
evaluate $\lambda_h$ specified in (\ref{eqn_labmdah}). Applying
(\ref{thm_eq_pe}) and (\ref{eq_pe_lnnn}) to (\ref{eqn_labmdah}), and
considering $\lim_{n \to \infty} \alpha_n = \alpha^* $ with
$|\alpha^{\star}| < \infty$, we establish
%\begin{subequations}
\begin{align}
 \lambda_h
 &  = n (h!)^{-1}(n q_n)^h
e^{-n q_n} \nonumber  \\
 &  \sim n (h!)^{-1} (\ln n)^h \cdot e^{-\ln n -
 (k-1)\ln \ln n - \alpha_n}  \nonumber  \\
 & = (h!)^{-1} (\ln n)^{h+1-k} e^{-\alpha_n}
 \nonumber  \\ & \to  \begin{cases} 0, &\textrm{ for }h = 0, 1,
\ldots,  k-2,  \\
 \frac{e^{-\alpha^*}}{(k-1)!},&\textrm{ for }h = k-1,
 \\ \infty ,&\textrm{ for }h = k, k+1, \ldots
\end{cases} \label{eqn_lbdh}
 \end{align}
By (\ref{eqn_phihell}) and (\ref{eqn_lbdh}), we derive that as $n
\to \infty$,
\begin{align}
\mathbb{P}[\phi_h = 0] & \to
\begin{cases} 1, & \textrm{ for }h = 0, 1,
\ldots,  k-2, \\
 e^{-\frac{e^{-\alpha^*}}{(k-1)!}}, & \textrm{ for }h = k-1,
 \\  0 , & \textrm{ for }h = k, k+1,
\ldots \end{cases}\label{eqn_expr_lahkk1}
 \end{align}
%
%Define $\lambda_h$ by $\lambda_h : = n (h!)^{-1}(n q_n)^h e^{-n
%q_n}$.  Then we use Lemma \ref{lem_pos_sum} in the following Charles
%Jordan's inversion formula \cite{JAZ:4922728}:
%\begin{align}
% &  \mathbb{P}[\phi_h = \ell] \nonumber  \\
% & \quad   = \Bigg\{ \sum_{m=\ell}^{n}(-1)^{m-\ell}
% \binom{m}{\ell} \times \nonumber  \\
% &   \times  \sum_{\begin{subarray}{c} (i_1, i_2, \ldots, i_m):
%   \\
%1 \leq i_1 < i_2 < \ldots < i_m \leq n
%\end{subarray}} \mathbb{P} [\textrm{Nodes }v_{i_1}, v_{i_2}, \ldots, v_{i_m}\textrm{ have degree
%}h] \Bigg\}, \nonumber
% \end{align}
%leading to
%\begin{align}
% \mathbb{P}[\phi_h = \ell]
% & \sim
% \sum_{m=\ell}^{\infty} \left[(-1)^{m-\ell} (m!)^{-1}
% \binom{m}{\ell} (\lambda_h)^m \right]
%\nonumber  \\
% & = \nonumber \sum_{m^*=0}^{\infty} \left[(-1)^{m^*}
%\frac{1}{\ell! m^*!}   (\lambda_h)^{\ell + m^*} \right] \\
% & = \frac{{\lambda_h} ^{\ell}e^{-\lambda_h}}{\ell !},
% \label{eqn_phihell}
% \end{align}
% Hence, $\phi_h$ asymptotically converges to a Poisson distribution
%with mean $\lambda_h$.

Using (\ref{eqn_expr_lahkk1}) in (\ref{eqn_1mindel2}) and
(\ref{eqn_1min}), we obtain $\mathbb{P}[\delta \geq k] \to
  e^{- \frac{e^{-\alpha ^*}}{(k-1)!}}$; i.e., Lemma
\ref{lem-mnd} is proved.

\subsection{Proving Lemma \ref{LEM1}} \label{prfseclem1}

We use $\mathcal {V}_m$ to denote the node set $\{v_1, v_2, \ldots,
v_m\}$. Lemma \ref{LEM1} evaluates the probability that each of
$\mathcal {V}_m$ has degree $h$. To compute such probability, we
look at whether at least two of $\mathcal {V}_{m}$ have an edge in
between, and whether at least two of $\mathcal {V}_{m}$ have at
least one common neighbor. To this end, we define $\mathcal {P}_1$ as the probability of event
\begin{align}
&  \begin{array}{c}
(\textrm{each of }\mathcal {V}_{m}
 \textrm{ has
degree }h)  \\ \bcap \big[  (\textrm{at least two of }\mathcal
{V}_{m}\textrm{ have an edge in between})\\ \bcup
(\textrm{at\hspace{2.5pt}least\hspace{2.5pt}two\hspace{2.5pt}of\hspace{2.5pt}}\mathcal
{V}_{m}\textrm{\hspace{2.5pt}have\hspace{2.5pt}at\hspace{2.5pt}least\hspace{2.5pt}one\hspace{2.5pt}common\hspace{2.5pt}neighbor})
\big],
\end{array}  \nonumber
 \end{align}
 and define $\mathcal {P}_2$ as the probability of event
 \begin{align}
& \begin{array}{c}
(\textrm{each of }\mathcal {V}_{m}
 \textrm{ has
degree }h)  \\ \bcap   (\textrm{no two of }\mathcal {V}_{m}\textrm{
have any edge in between})\\ \bcap   (\textrm{no two of }\mathcal
{V}_{m}\textrm{ have any common neighbor}).
\end{array} \nonumber
 \end{align}
Then $\mathbb{P} [\textrm{each of }\mathcal {V}_{m}
 \textrm{ has
degree }h] = \mathcal {P}_1 + \mathcal {P}_2$. Thus, Lemma
\ref{LEM1} will hold once we establish the following two
propositions.

\begin{proposition}\label{PROP_ONE}
Under the conditions of Theorem \ref{THM1},
it holds that $\mathcal {P}_1=  o \left((h!)^{-m} (n q_n)^{hm} e^{-m n
q_n}\right)$.
%\begin{align}
%& \mathbb{P} [\mathcal {E}_1]=  o \left((h!)^{-m} (n q_n)^{hm} e^{-m n
%q_n}\right).  \nonumber
% \end{align}

\end{proposition}

\begin{proposition} \label{PROP_SND}
Under the conditions of Theorem \ref{THM1},
it holds that $\mathcal {P}_2 \sim  (h!)^{-m}  (n q_n)^{hm} e^{-m n q_n}$.
%\begin{align}
%& \mathbb{P} [\mathcal {E}_2] \sim  (h!)^{-m}  (n q_n)^{hm} e^{-m n q_n}.  \nonumber
% \end{align}
\end{proposition}

To prove Propositions \ref{PROP_ONE} and \ref{PROP_SND}, we analyze
below how nodes in graph $\mathbb{G} $ have edges.
We first look at how edges exist between $v_1, v_2, \ldots, v_m$.
Recalling ${C}_{i j}$ as the event that the communication channel
between distinct nodes $v_i$ and $v_j$ is {\em on}, we set
$\boldsymbol{1}[C_{ij}]$ as the indicator variable of event ${C}_{i
j}$ by
\begin{align}
 \hspace{-2pt} \boldsymbol{1}[C_{ij}]& \hspace{-2pt} :=  \hspace{-2pt} \begin{cases}
1,~ \textrm{if the
channel between }v_i\textrm{ and }v_j\textrm{ is \textit{on}}, \\
0,~ \textrm{if the channel between }v_i\textrm{ and }v_j\textrm{ is
\textit{off}}.
\end{cases} \nonumber
\end{align}
We denote by $\mathcal {C}_m$ a $\binom{m}{2}$-tuple consisting of
all possible $\boldsymbol{1}[C_{ij}]$ with $1 \leq i < j \leq m$ as
follows:\vspace{-2pt} 
\begin{align}
\mathcal {C}_m : = ( &\boldsymbol{1}[C_{12}],
,\ldots,\boldsymbol{1}[C_{1m}],~~~\boldsymbol{1}[C_{23}],
,\ldots,\boldsymbol{1}[C_{2m}], \nonumber \\ &
  \boldsymbol{1}[C_{34}], \ldots,\boldsymbol{1}[C_{3m}],~~~\ldots,~~~
\vspace{-2pt} \boldsymbol{1}[C_{(m-1),m}]). \nonumber
\end{align}
Recalling $S_i$ as the key set on node $v_i$, we define a $m$-tuple
$\mathcal {T}_m$ through $\mathcal {T}_m  : = (S_1, S_2, \ldots,
S_m)$. Then we define $\mathcal {L}_m$ as $\mathcal {L}_m   : =
(\mathcal {C}_m, \mathcal {T}_m)$. With $\mathcal {L}_m$, we have
the \emph{on}/\emph{off} states of all channels between nodes $v_1,
v_2, \ldots, v_m$ and the key sets $S_1, S_2, \ldots, S_m$ on these
$m$ nodes, so all edges between these $m$ nodes in graph
$\mathbb{G}\iffalse_{on}\fi$ are determined. Let $\mathbb{C}_m,
\mathbb{T}_m$ and $\mathbb{L}_m$ be the sets of all possible
$\mathcal {C}_m, \mathcal {T}_m$ and $\mathcal {L}_m$, respectively.

Now we further introduce some notation to characterize how nodes
$v_1, v_2, \ldots, v_m$ have edges with nodes of $\overline{\mathcal
{V}_m}$, where $\overline{\mathcal {V}_m}$ denotes $\{v_{m+1},
v_{m+2}, \ldots, v_n\}$. Let $N_i$ be the neighborhood set of node
$v_i$, i.e., the set of nodes that have edges with $v_i$.
%
%We define  Recalling $\mathcal {V}_m = \{v_1, v_2, \ldots, v_m\}$,
%we denote $ \{v_1, v_2, \ldots, v_n\} \setminus \mathcal {V}_m $ by
%$\overline{\mathcal {V}_m}$.
 We also define set $\overline{N_i}$ as the set $\{v_{m+1}, v_{m+2},
\ldots, v_n\} \setminus N_i$. Then we are ready to define sets
$M_{j_1 j_2 \ldots j_m}$
 for all $j_1, j_2, \ldots, j_m \in \{0,1\}$ which characterize the
 relationships between sets $N_i$ for
$i=1,2,\ldots,m$. We define\vspace{-2pt} 
\begin{align}
&  M_{j_1 j_2 \ldots j_m} \hspace{-1pt} : \vspace{-2pt} = \hspace{-1pt}
\bigg(\hspace{-1pt} \bigcap_{i\in \{1,2,\ldots,m\}:j_{i}=1} N_{i}
\hspace{-1pt} \bigg)\hspace{-1pt} \bcap \hspace{-1pt}
\bigg(\hspace{-1pt} \bigcap_{i\in \{1,2,\ldots,m\}:j_{i}=0}
\overline{N_{i}} \hspace{-1pt} \bigg). \label{Mjdefi}
\end{align}
In other words, for $i=1,2,\ldots,m$, if $N_i$ is not empty, each
node in $N_i$ belongs to $M_{j_1 j_2 \ldots j_m}$ if $j_i = 1$ and
does not belong to $M_{j_1 j_2 \ldots j_m}$ if $j_i = 0$. Also, if
$j_1 = j_2 = \ldots = j_m = 0$, then $M_{j_1 j_2 \ldots j_m} =
\bigcap_{i=1}^{m}\overline{N_{i}}$. The sets $M_{j_1 j_2 \ldots
j_m}$ for $j_1, j_2, \ldots, j_m \in \{0,1\}$ are mutually disjoint,
and constitute a partition of the set $\overline{\mathcal {V}_m}$ (a
partition is allowed to contain empty sets here). By the definition
of $M_{j_1 j_2 \ldots j_m}$ for $j_1, j_2, \ldots, j_m \in \{0,1\}$,
we \vspace{-2pt} have
\begin{align}
\sum_{j_1, j_2, \ldots, j_m \in \{0,1\}}|M_{j_1 j_2 \ldots j_m}^{*}|
= |\overline{\mathcal {V}_m}| = \vspace{-2pt} n - m, \label{eqn_nodevi_h5}
\end{align}
and\vspace{-2pt} 
\begin{align}
\hspace{-5pt} \sum_{\begin{subarray}{c}j_1, j_2, \ldots, j_m \in \{0,1\}: \\
\sum_{i=1}^{m}j_i \geq 1.
\end{subarray}}\hspace{-2pt} |M_{j_1 j_2 \ldots j_m} | =
\bigg|\bigg( \bigcup_{i=1}^m N_i \bigg) \bcap \overline{\mathcal
{V}_m}\bigg|\vspace{-2pt}  .\label{eqn_nodevi_h2}
\end{align}
%We will use (\ref{eqn_nodevi_h5}) and (\ref{eqn_nodevi_h2}) later.

We further define $2^m$-tuple\vspace{-2pt}  $\mathcal {M}_m$ through\footnote{For
a non-negative integer $x$, the term $0^{x}$ is short for
$\underbrace{00 \ldots 0}_{\textrm{``}x\textrm{''} \textrm{ number
of ``}0\textrm{''}}$. Also, for clarity, we add commas in the
subscript of $M_{0^{m-2}1,0}$ etc.}
%\footnote{We use $0^{x}$ as an abbreviation of $\underbrace{00
%\ldots 0}_{\textrm{``}x\textrm{''} \textrm{ number of
%``}0\textrm{''}}$ with a non-negative integer $x$.}
\begin{align}
 \mathcal {M}_m & = \vspace{-2pt}   \big( | M_{j_1 j_2 \ldots j_m} |
 \hspace{2pt}\boldsymbol{\big|}  \hspace{2.5pt} j_1, j_2, \ldots, j_m \in \{0,1\} \big)  \nonumber  \\
&   = \vspace{-2pt}   \big( | M_{0^m} |, |M_{0^{m-1},1}|, |M_{0^{m-2}1,0}| ,
|M_{0^{m-2}1,1}| , \ldots\big),\nonumber
\end{align}
where $| M_{j_1 j_2 \ldots j_m} |$ means the cardinality of $M_{j_1 j_2 \ldots j_m}$.
Under event $\mathcal {E}_2$, the set $\mathcal {M}_m$ is determined
and we denote its value by $\mathcal {M}_m^{(0)}$, which satisfies\vspace{-2pt} 
\begin{align}
\begin{cases}
|M_{0^{i-1}, 1, 0^{m-i}}|   = h,\vspace{-2pt} & \textrm{for }i=1,2,\ldots,m; \\
|M_{j_1 j_2 \ldots j_m}|  = 0,\vspace{-2pt} &\textrm{for } \sum_{i=1}^m j_i >
1;\\|M_{0^m}|
  = n - m - hm .\vspace{-2pt} 
\end{cases}\label{mmsizesat} \vspace{-2pt}
\end{align}
To analyze event $\mathcal {E}_2$, we define $\mathbb {L}_m^{(0)}$ such that $\big(\mathcal
{L}_m \in \mathbb {L}_m^{(0)}\big)$ is the event that no two of
nodes $ v_1, v_2, \ldots, v_m $ have any common neighbor. In view of events $\big( \mathcal {L}_m \in \mathbb{L}_m^{(0)}
\big)$, $\big( \mathcal {M}_m = \mathcal{M}_m^{(0)} \big)$ and
$\mathcal {E}_2$, then $\mathcal {E}_2$ is the same as
$\big( \mathcal {L}_m \in \mathbb{L}_m^{(0)} \big) \bcap \big(
\mathcal {M}_m = \mathcal{M}_m^{(0)} \big)$; \vspace{-2pt}  i.e.,
\begin{align}
\mathcal {E}_2 & = \big[\big( \mathcal {L}_m \in \mathbb{L}_m^{(0)}
\big) \bcap \big( \mathcal {M}_m \vspace{-2pt} = \mathcal{M}_m^{(0)} \big) \big].
\label{e2}
\end{align}

We define $\mathbb{M}_m(\mathcal {L}_m)$ for $\mathcal {L}_m \in
\mathbb {L}_m $ as the set of $\mathcal {M}_m$ under which each of
$\mathcal {V}_{m}$ has degree $h$. Thus, the event that each of
$\mathcal {V}_{m}$ has degree $h$ is $\big(\mathcal {L}_m \in
\mathbb {L}_m \big) \bcap \big( \mathcal {M}_m \in
\mathbb{M}_m(\mathcal {L}_m) \big)$, which together with (\ref{e2})
yields \vspace{-2pt} 
\begin{align}
\mathcal {E}_1 \hspace{-2pt} \vspace{-100pt}  =  \hspace{-2pt} \bigcup_{\begin{subarray}{c}
\mathcal {L}_m^{*} \in \mathbb{L}_m,
\hspace{2pt} \mathcal {M}_m^{*} \in \mathbb{M}_m (\mathcal {L}_m^{*}):  \\
   \left(\mathcal {L}_m^{*}
   \notin \mathbb{L}_m^{(0)}\right)
\textrm{ or }\left(\mathcal {M}_m^{*} \neq
\mathcal{M}_m^{(0)}\right)
\end{subarray}} \hspace{-2pt}
 \mathbb{P} \big[ \big( \mathcal {L}_m \hspace{-1pt} =
\hspace{-1pt} \mathcal {L}_m^{*} \big) \hspace{-2pt} \bcap \hspace{-2pt} \big( \mathcal {M}_m
\hspace{-1pt} = \hspace{-1pt} \mathcal {M}_m^{*} \big) \big].\vspace{-10pt} 
\label{e1}
\end{align}

Now we prove Propositions \ref{PROP_ONE} and \ref{PROP_SND} based on
(\ref{e2}) and (\ref{e1}). The inequality below following from (\ref{eq_pe_lnnn}) will be applied\vspace{-1pt}  often:\vspace{-2pt} 
\begin{align}
q_n & \leq \frac{2\ln n}{n}\textrm{ for all $n$ sufficiently large}.\vspace{-2pt} 
\label{eq_pe_upper}
\end{align}
\iffalse

Also, we will   use
 $1+x \leq e^x$ for any
real $x$ and \cite[Fact 2]{ZhaoYaganGligor} which says $1 - xy \leq (1-x)^y \leq 1 - xy +
\frac{1}{2} x^2 y ^2$ for $0\leq x <1$ and $y = 0 , 1, 2, \ldots$.\vspace{-2pt} 

aaa

\fi

%\begin{align}
% & \mathcal {M}_m^{(0)} \nonumber  \\
%& \quad  =\big\{ \mathcal {M}_m \boldsymbol{\mid}|M_{0^{i-1}, 1,
%0^{m-i}}| = h \textrm{ for }i=1,2,\ldots,m, \nonumber  \\
%& \quad \quad \quad \quad \quad \hspace{7pt} |M_{j_1 j_2 \ldots
%j_m}| = 0\textrm{ for any }\sum_{i=1}^m
%j_i > 1,  \nonumber  \\
%& \quad \quad \quad \quad \quad \hspace{7pt}  M_{0^m}
% = n - m - hm . \big\} . \nonumber
%\end{align}

\subsubsection{The Proof of Proposition \ref{PROP_ONE}}~
\label{sec:PROP_ONE}\vspace{-2pt} 

In view of (\ref{e1}) and considering the disjointness of events
$\big( \mathcal {L}_m   =   \mathcal {L}_m^{*} \big) \bcap \big(
\mathcal {M}_m  =   \mathcal {M}_m^{*} \big)$ for $\mathcal
{L}_m^{*} \in \mathbb{L}_m$ and $\mathcal {M}_m^{*} \in \mathbb{M}_m
(\mathcal {L}_m^{*})$, we express $\mathbb{P}
 [\mathcal {E}_1 ]$ as\vspace{-2pt} 
  \begin{align}
 \sum_{\begin{subarray}{c}
\mathcal {L}_m^{*} \in \mathbb{L}_m,
\hspace{2pt} \mathcal {M}_m^{*} \in \mathbb{M}_m (\mathcal {L}_m^{*}):  \\
   \left(\mathcal {L}_m^{*}
   \notin \mathbb{L}_m^{(0)}\right)
\textrm{ or }\left(\mathcal {M}_m^{*} \neq
\mathcal{M}_m^{(0)}\right)
\end{subarray}}  \hspace{-2pt}
 \mathbb{P} \big[ \big( \mathcal {L}_m \hspace{-2pt} =
\hspace{-2pt} \mathcal {L}_m^{*} \big) \hspace{-2pt} \bcap \hspace{-2pt} \big( \mathcal {M}_m
\hspace{-2pt} = \hspace{-2pt} \mathcal {M}_m^{*} \big)
\big]\label{term1}\vspace{-2pt}
\end{align}

We evaluate (\ref{term1}) by computing\vspace{-2pt} 
\begin{align}
\mathbb{P} \big[ \big( \mathcal {M}_m = \mathcal {M}_m^{*} \big)
\boldsymbol{\mid} \mathcal {L}_m = \mathcal {L}_m^{*} \big].\vspace{-2pt} 
\label{eq_MmMm}
\end{align}
With $\mathcal {C}_m ^{*}$ and $\mathcal {T}_m ^{*}$ defined such
that $\mathcal {L}_m^{*} = (\mathcal {C}_m^{*}, \mathcal
{T}_m^{*})$, event $(\mathcal {L}_m \hspace{-1pt}=\hspace{-1pt} \mathcal {L}_m^{*})$
is the union of events $(\mathcal {C}_m \hspace{-1pt}=\hspace{-1pt} \mathcal
{C}_m^{*})$ and $(\mathcal {T}_m \hspace{-1pt}=\hspace{-1pt} \mathcal {T}_m^{*})$.\\
Since $( \mathcal {C}_m \hspace{-2pt} = \hspace{-2pt} \mathcal
{C}_m^{*} )$ and $( \mathcal {M}_m \hspace{-2pt} = \hspace{-2pt}
\mathcal {M}_m^{*} )$ are independent, we\vspace{-2pt}  get
\begin{align}
(\ref{eq_MmMm}) & = \mathbb{P} \big[ \big( \mathcal {M}_m = \mathcal
{M}_m^{*} \big)
 \boldsymbol{\mid}
\big( \mathcal {T}_m = \mathcal {T}_m^{*} \big) \big] \nonumber
.\vspace{-2pt}
\end{align}

For each $j_1, j_2, \ldots, j_m \in \{0,1\},$ for any distinct nodes
$w_1 , w_2  \in  \overline{\mathcal {V}_m} $,
events $(w_1 \hspace{-2pt} \in \hspace{-2pt} M_{j_1 j_2 \ldots
j_m})$ and $(w_2 \in M_{j_1 j_2 \ldots j_m})$ are  conditionally independent given $(\mathcal
{T}_m = \mathcal {T}_m^{*})$ , where $\mathcal {T}_m^{*}$
specifies the key sets $S_1, S_2, \ldots,
S_m$ as $S_1^{*}, S_2^{*}, \ldots, S_m^{*}$,
respectively). Thus, with $\mathcal {M}_m^{*}  $ being $   \big( |
M_{0^m}^{*} |, |M_{0^{m-1},1}^{*}|, |M_{0^{m-2}1,0}^{*}| ,
|M_{0^{m-2}1,1}^{*}| , \ldots\big)$, we obtain
\begin{align}
 & \hspace{-2pt} (\ref{eq_MmMm}) = f(n-m , \mathcal {M}_m^{*})\mathbb{P} [w \in M_{0^m}
\hspace{-2pt} \boldsymbol{\mid} \hspace{-2pt} \mathcal
{T}_m = \mathcal {T}_m^{*} ]^{|M_{0^m}^{*} |} \times \vspace{-2pt} \nonumber  \\
& \hspace{-2pt} \vspace{-10pt} \quad \prod_{\begin{subarray}{c}j_1, j_2, \ldots, j_m \in \{0,1\}: \\
\sum_{i=1}^{m}j_i \geq 1.
\end{subarray}} \hspace{-2pt} \mathbb{P}[w \in M_{j_1 j_2 \ldots j_m}
 \hspace{-2pt} \boldsymbol{\mid} \hspace{-2pt} \mathcal
{T}_m = \mathcal {T}_m^{*}]^{|M_{j_1 j_2 \ldots j_m}^{*} |},\vspace{-5pt} 
\label{eqn_epsilonmmmmold}
\end{align}
where $f(n-m , \mathcal {M}_m^{*})$ is the number of ways assigning the $(n-m)$ nodes from $ \overline{\mathcal {V}_m} $ to $M_{j_1 j_2 \ldots j_m}$ such that  $|M_{j_1 j_2 \ldots j_m}|$ equals $|M_{j_1 j_2 \ldots j_m}^{*} |$, for $j_1, j_2, \ldots, j_m \in \{0,1\}$. Then\vspace{-2pt} 
\begin{align}
&  f(n-m  , \mathcal {M}_m^{*}) = \frac{ (n-m ) !}{\prod_{j_1, j_2,
\ldots, j_m \in \{0,1\}}
(|M_{j_1 j_2 \ldots j_m}^{*}|! )},\vspace{-1pt}  \label{eqn_fnexpr}
\end{align}
which along with (\ref{eqn_nodevi_h5}) yields\vspace{-1pt} 
\begin{align}
f(n-m  , \mathcal {M}_m^{*})    \hspace{-1pt} &   \leq \hspace{-1pt}   [(n-m ) !] / (|M_{0^m}^{*} |!) \vspace{-2pt} \nonumber %  \\
%&  \hspace{-2pt} = \hspace{-2pt}   (n \hspace{-2pt}
%-\hspace{-2pt} m ) ! \hspace{-2pt} \Big/ \hspace{-2pt}
%\Big(n\hspace{-2pt} -\hspace{-2pt} m - \hspace{-3pt}
%\sum_{\begin{subarray}{c}j_1, j_2, \ldots, j_m \in \{0,1\}: \\
%\sum_{i=1}^{m}j_i \geq 1.
%\end{subarray}}|M_{j_1 j_2 \ldots j_m}^{*}|\hspace{-1pt}\Big)!
 \\
 \hspace{-1pt} &   \leq \hspace{-1pt} \textstyle{ n^{\sum_{\begin{subarray}{c}j_1,
j_2, \ldots, j_m
\in \{0,1\}: \\
\sum_{i=1}^{m}j_i \geq 1.
\end{subarray}}|M_{j_1 j_2 \ldots j_m}^{*}|}.}\vspace{-2pt}  \label{eqn_fnm}
\end{align}

For any $j_1, j_2, \ldots, j_m \in \{0,1\}$ with $\sum_{i=1}^{m}j_i
\geq 1$, there exists $t \in \{0,1,\ldots, m\}$ such that $j_t = 1$,
so \vspace{-2pt}
\begin{align}
&  \mathbb{P}\big[w \in M_{j_1 j_2 \ldots j_m} \boldsymbol{\mid}
\mathcal {T}_m = \mathcal {T}_m^{*} \big] \vspace{-2pt} \nonumber  \\
&  \quad \leq  \mathbb{P}[E_{w v_t}  \boldsymbol{\mid} \mathcal {T}_m =
\mathcal {T}_m^{*}] = \mathbb{P}[E_{w v_t} ] = q_n, \vspace{-2pt}
\label{eqn_pe_not00}
\end{align}
where $E_{w v_t}$ is the event that  an edge exists between
nodes $w$ and $v_t$. Substituting (\ref{eqn_fnm}) and (\ref{eqn_pe_not00}) into
(\ref{eqn_epsilonmmmmold}), and denoting $\sum_{\begin{subarray}{c}j_1, j_2, \ldots, j_m \in \{0,1\}: \\
\sum_{i=1}^{m}j_i \geq 1.
\end{subarray}}|M_{j_1 j_2 \ldots j_m}^{*}|$ by $\Lambda$, we obtain
\begin{align}
\hspace{-1pt}  (\ref{eq_MmMm})
 & <
 (nq_n)^{\Lambda}
\hspace{-1pt} \times \hspace{-1pt} \mathbb{P} [w \in M_{0^m}
\hspace{-1pt} \boldsymbol{\mid}  \hspace{-1pt} \mathcal {T}_m =
\mathcal {T}_m^{*} ]^{| M_{0^m}^{*} |  } . \label{eqn_epsilonmmmm} \vspace{-2pt}
\end{align}

To further evaluate (\ref{eq_MmMm}) based on (\ref{eqn_epsilonmmmm}),
we will prove below that if $\big(\mathcal {L}_m^{*} \notin
\mathbb{L}_m^{(0)} \big)$ or $\big(\mathcal {M}_m^{*} \neq
\mathcal{M}_m^{(0)} \big)$, then
\begin{align}
\Lambda \leq hm - 1 . \vspace{-2pt} \label{lambda}
\end{align}

 On the one hand, if $\mathcal
{L}_m^{*} \notin \mathbb{L}_m^{(0)}$, there exist $i_1$ and $i_2$
with $1 \leq i_1 < i_2 \leq m$ such that nodes $v_{i_1}$ and
$v_{i_2}$ are neighbors. Hence, $ \{v_{i_1},
v_{i_2}\} \subseteq [( \bigcup_{i=1}^m N_i ) \bigcap \mathcal {V}_m
]$ holds. Then from (\ref{eqn_nodevi_h2}), we have $\Lambda =
 \big|\bigcup_{i=1}^m N_i\big|  -
  \big|\big( \bigcup_{i=1}^{m} N_i \big) \bcap \mathcal {V}_m\big|
 \leq hm - 2.$
%\begin{align}
% & \Lambda =
% \big|\bigcup_{i=1}^m N_i\big|  -
%  \big|\big( \bigcup_{i=1}^{m} N_i \big) \bcap \mathcal {V}_m\big|
% \leq hm - 2. \nonumber %\label{eq_first}
%\end{align}
 On the other hand, if $\mathcal {M}_m^{*} \neq \mathcal{M}_m^{(0)}$,
there exist $i_3 $ and $ i_4$ with $1 \leq i_3 < i_4 \leq m$ such
that $N_{i_3} \cap N_{i_4} \neq \emptyset$. Then from
(\ref{eqn_nodevi_h2}), $\Lambda \leq
 \big|\bigcup_{i=1}^m N_i\big|  \leq
  \big(\sum_{i=1}^m |N_i|\big) - |N_{i_3} \cap N_{i_4}| \leq
 hm - 1$ follows. Thus, we have proved (\ref{lambda}),
%\begin{align}
% & \Lambda \leq
% \big|\bigcup_{i=1}^m N_i\big|  \leq
%  \big(\sum_{i=1}^m |N_i|\big) - |N_{i_3} \cap N_{i_4}| \leq
% hm - 1. \nonumber %\label{eq_first}
%\end{align}
which along with (\ref{eqn_nodevi_h5}) \vspace{-2pt}  leads to
\begin{align}
| M_{0^m}^{*} |  &=  n - m -
 \Lambda > n - m - hm. \vspace{-2pt} \label{eqn_M0m2n}
\end{align}

From (\ref{eq_pe_lnnn}), it is true that $nq_n \sim \ln n$, implying $nq_n > 1$ for
all $n$ sufficiently large. Then substituting (\ref{lambda}) and (\ref{eqn_M0m2n}) into
(\ref{eqn_epsilonmmmm}), we obtain that if $\big(\mathcal
{L}_m^{*} \notin \mathbb{L}_m^{(0)} \big)$ or $\big(\mathcal
{M}_m^{*} \neq \mathcal{M}_m^{(0)} \big)$, then for
all $n$ sufficiently large, it holds that\vspace{-2pt} 
%(\ref{eqn_fnm}) (\ref{lambda}) (\ref{eqn_M0m2n}) and
%(\ref{eqn_pe_not00})
\begin{align}
\hspace{-1pt}  (\ref{eq_MmMm})
 & <
 (nq_n)^{hm - 1}
\hspace{-1pt} \times \hspace{-1pt} \mathbb{P} [w \in M_{0^m}
\hspace{-1pt} \boldsymbol{\mid}  \hspace{-1pt} \mathcal {T}_m =
\mathcal {T}_m^{*} ]^{n - m - hm} . \vspace{-2pt}  \label{eq_pmmll2}
\end{align}

Applying (\ref{eq_MmMm}) and (\ref{eq_pmmll2}) to (\ref{term1}),% and
%considering $\mathbb{P} \big[ \mathcal {L}_m = \mathcal {L}_m^{*}
%\big] = \mathbb{P} \big[ \mathcal {T}_m = \mathcal {T}_m^{*} \big]
%\mathbb{P} \big[ \mathcal {C}_m = \mathcal {C}_m^{*} \big]$,
 we
get\vspace{-2pt} 
\begin{align}
(\ref{term1}) & \hspace{-1pt} \vspace{-2pt} < \hspace{-1pt} \sum_{\mathcal {L}_m^{*} \in \mathbb{L}_m} \hspace{-2pt}
\Big\{ |\mathbb{M}_m (\mathcal {L}_m^{*})| \hspace{-1pt} \times \hspace{-1pt} \mathbb{P} \big[ \mathcal {L}_m \hspace{-1pt} = \hspace{-1pt} \mathcal {L}_m^{*} \big]  \hspace{-1pt} \times \hspace{-1pt} \textrm{R.H.S. of (\ref{eq_pmmll2})} \Big\} .\label{eqn_TmCmt}
\end{align}

%\begin{align}
%&  \mathbb{P} \big[ \big( \mathcal {M}_m = \mathcal {M}_m^{*} \big)
%\boldsymbol{\mid} \mathcal {L}_m = \mathcal {L}_m^{*} \big]
%\nonumber  \\
%& = \sum_{\mathcal {T}_m^{*} \in \mathbb{T}_m} \mathbb{P} \big[
%\big( \mathcal {M}_m = \mathcal {M}_m^{*} \big) \cap \big( \mathcal
%{T}_m = \mathcal {T}_m^{*} \big) \boldsymbol{\mid} \mathcal {L}_m =
%\mathcal {L}_m^{*} \big] \nonumber  \\
%& \leq 2 (nq_n)^{\sum_{\begin{subarray}{c}j_1, j_2, \ldots, j_m \in \{0,1\}, \\
%\sum_{i=1}^{m}j_i \geq 1.
%\end{subarray}}|M_{j_1 j_2 \ldots j_m}^{*}|}
% e^{-m n q_n}   . \label{eq_pmmll}
%\end{align}

%
%
%Then
%\begin{align}
%(\ref{eq_except00}) &   \leq \sum_{\mathcal {L}_m^{*} \in
%\mathbb{L}_m}  2(nq_n)^{hm - 1}
% e^{-m n q_n} |\mathbb{M}_m (\mathcal {L}_m^{*})| \mathbb{P}
%\big[ \mathcal {L}_m = \mathcal {L}_m^{*} \big]
%\label{eq_except00_boundmm}
%\end{align}
To bound $|\mathbb{M}_m (\mathcal {L}_m^{*})|$, note that $\mathcal
{M}_m$ is a $2^m$-tuple. Among the $ 2^m $ elements of the tuple,
each of $|M_{j_1 j_2 \ldots j_m}
|\big|_{\begin{subarray}{c}j_1, j_2, \ldots, j_m \in \{0,1\}: \\
\sum_{i=1}^{m}j_i \geq 1.
\end{subarray}}$ is at least 0 and at most $h$; and the remaining
element $| M_{0^m} |$ can be determined by (\ref{eqn_nodevi_h5}).
Then it's straightforward that $|\mathbb{M}_m (\mathcal {L}_m^{*})|   \leq (h+1)^{2^m-1}$.
%\begin{align}
%|\mathbb{M}_m (\mathcal {L}_m^{*})|   \leq (h+1)^{2^m-1}.
%\label{eqn_MmLm}
%\end{align}
Using this result in (\ref{eqn_TmCmt}), and considering
$\big(\mathcal {L}_m = \mathcal {L}_m^{*}\big)$ is the union of
independent events $\big(\mathcal {T}_m = \mathcal {T}_m^{*}\big)$
and $\big(\mathcal {C}_m \hspace{-1pt} = \hspace{-1pt} \mathcal
{C}_m^{*}\big)$, and $\sum_{\mathcal {C}_m^{*} \in \mathbb{C}_m}
\hspace{-1pt} \mathbb{P} \big[ \mathcal {C}_m \hspace{-1pt} =
\hspace{-1pt} \mathcal {C}_m^{*} \big] \hspace{-1pt} = \hspace{-1pt}
1$, we derive\vspace{-2pt} 
\begin{align}
(\ref{term1}) &  \vspace{-2pt} <  (h+1)^{2^m-1} (nq_n)^{hm-1} \hspace{-2pt} \times
\hspace{-2pt}\sum_{\mathcal {T}_m^{*} \in \mathbb{T}_m}
\hspace{-4pt} \Big\{
\mathbb{P}\big[ \mathcal {T}_m = \mathcal {T}_m^{*} \big] \vspace{-2pt}  \nonumber  \\
& \vspace{-2pt}  \quad \times \mathbb{P} [w \in M_{0^m} \boldsymbol{\mid} \mathcal
{T}_m = \mathcal {T}_m^{*} ]^{n - m - hm} \Big\} . \vspace{-2pt}
\label{prop_prf}
\end{align}
From (\ref{prop_prf}) and $nq_n \sim \ln n \to \infty$ as $n\to\infty$ by
(\ref{eq_pe_lnnn}), the proof of Proposition \ref{PROP_ONE} is
completed once we show \vspace{-2pt}
\begin{align}
 & \sum_{ \mathcal {T}_m^{*} \in \mathbb{T}_m } \mathbb{P}[\mathcal {T}_m = \mathcal {T}_m^{*}]
\mathbb{P} [w \in M_{0^m} \boldsymbol{\mid} \mathcal
{T}_m = \mathcal {T}_m^{*} ]^{n - m - hm} \vspace{-2pt} \nonumber  \\
& \quad  \leq e^{- m n q_n} \cdot [1+o(1)] . \vspace{-2pt}   \label{EQ}
\end{align}

\subsection{Establishing (\ref{EQ})}\vspace{-2pt}

From (\ref{eq_evalprob_3}) and (\ref{eq_evalprob_4}) (Lemma \ref{lem_evalprob}
in the Appendix),  we get\vspace{-2pt} 
\begin{align}
 & \mathbb{P} [w \in M_{0^m}^{*} \boldsymbol{\mid} \mathcal
{T}_m = \mathcal {T}_m^{*} ]^{n - m - hm } \vspace{-2pt} \nonumber  \\
& \hspace{-5pt}  \vspace{-2pt} = \hspace{-2pt} \mathbb{P} [w \hspace{-2pt} \in
\hspace{-2pt} M_{0^m}^{*} \hspace{-2pt} \boldsymbol{\mid}
\hspace{-2pt} \mathcal {T}_m \hspace{-2pt} = \hspace{-2pt} \mathcal
{T}_m^{*} ]^{ n} \mathbb{P} [w \hspace{-2pt} \in \hspace{-2pt}
M_{0^m}^{*} \hspace{-2pt} \boldsymbol{\mid}
\hspace{-2pt} \mathcal {T}_m \hspace{-2pt} = \hspace{-2pt} \mathcal {T}_m^{*} ]^{-m - h m}  \nonumber  \\
& \hspace{-5pt} \leq \vspace{-2pt}  \hspace{-2pt} e^{- m n q_n \hspace{-.5pt} +
\hspace{-.4pt} m^2 n {q_n}^2 \hspace{-2pt} + \hspace{-1pt}
  \frac{n q_n p_n}{K_n}\hspace{-2pt}\sum_{1\leq i <j \leq m}\hspace{-2pt}|S_{ij}^{*}|}
   \hspace{-1pt} (1\hspace{-1.5pt}-\hspace{-1.5pt}m q_n)^{-m - h m}\vspace{-2pt} \label{eqn_prbwM}
\end{align}
for all $n$ sufficiently large,
where $S_{ij}^{*} := S_{i}^{*} \cap S_{j}^{*}$. With
(\ref{eq_pe_lnnn}) (i.e., $q_n \sim \frac{\ln n}{n}$), we have $m^2
n {q_n}^2
  = o(1)$ and $m q_n =
o(1)$, which are substituted into (\ref{eqn_prbwM}) to induce
(\ref{EQ}) once we prove\vspace{-2pt} %
%
%
%obtain
%\begin{align}
% & \sum_{ \mathcal {T}_m^{*} \in \mathbb{T}_m  } \mathbb{P}[\mathcal {T}_m = \mathcal {T}_m^{*}]
%\mathbb{P} [w \in M_{00 \ldots 0}^{*} \boldsymbol{\mid} \mathcal
%{T}_m = \mathcal {T}_m^{*} ]^{n - m - hm} \nonumber  \\
%& \quad  \leq e^{- m n q_n} \cdot [1+o(1)] \nonumber  \\
%& \quad\quad \times \sum_{ \mathcal {T}_m^{*} \in \mathbb{T}_m  }
%\mathbb{P}[\mathcal {T}_m = \mathcal {T}_m^{*}] e^{\frac{n q_n
%p_n}{K_n}\sum_{1\leq i <j \leq m}|S_{ij}|}. \label{eqn_sum_Tm}
%\end{align}
%By (\ref{eqn_sum_Tm}), establishing (\ref{EQ}) is equivalent to
%proving
\begin{align}
\sum_{ \mathcal {T}_m^{*} \in \mathbb{T}_m  } \mathbb{P}[\mathcal
{T}_m = \mathcal {T}_m^{*}] e^{\frac{n q_n p_n}{K_n}\sum_{1\leq i <j
\leq m}|S_{ij}^{*}|} &  \leq 1+o(1).\vspace{-2pt}  \label{eqn_sumTmst}
\end{align}
 L.H.S. of (\ref{eqn_sumTmst}) is denoted by $H_{n,m}$ and evaluated
below. For each fixed and sufficiently large $n$, we consider: {a)}
{${ p_n <  n^{-\delta} (\ln n)^{-1}}$} and {b)} {${ p_n \geq
n^{-\delta} (\ln n)^{-1}}$}, where $\delta$ is an arbitrary constant
with $0<\delta<1$.

\noindent \textbf{a)} $\boldsymbol{ p_n <  n^{-\delta} (\ln
n)^{-1}}$

From $p_n < n^{-\delta} (\ln n)^{-1}$, $|S_{ij}^{*}| \leq K_n$ for $1\leq
i <j \leq m$ and (\ref{eq_pe_upper}), then for all $n$ sufficiently large, it holds that
$e^{\frac{n q_n p_n}{K_n} \sum_{1\leq i <j
\leq m}|S_{ij}^{*}|}  \hspace{-1.5pt}<\hspace{-1.5pt} e^{2 n^{-\delta}\hspace{-1pt}  \cdot \hspace{-1pt}\binom{m}{2}} \hspace{-1.5pt}< \hspace{-1.5pt} e^{ m^2
n^{-\delta}}$,
which is used in $H_{n,m}$ so that $H_{n,m} \hspace{-1.5pt}<\hspace{-1.5pt} e^{ m^2 n^{-\delta}} \hspace{-1pt}\sum_{ \mathcal {T}_m^{*} \in
\mathbb{T}_m  } \mathbb{P}[\mathcal {T}_m \hspace{-1pt}=\hspace{-1pt} \mathcal {T}_m^{*}] \hspace{-1.5pt}=\hspace{-1.5pt}
e^{ m^2 n^{-\delta}} $.
%\begin{align}
%& H_{n,m} < e^{ m^2 n^{-\delta}} \sum_{ \mathcal {T}_m^{*} \in
%\mathbb{T}_m  } \mathbb{P}[\mathcal {T}_m = \mathcal {T}_m^{*}] =
%e^{ m^2 n^{-\delta}} .\nonumber
%\end{align}

% $F_{m-1}$ is the event that $S_{1}, S_{2}, \ldots,
%S_{m-1}$ are mutually disjoint, meaning that $T_{j_1 j_2 \ldots
%j_{m-1}} = \emptyset$ for any $ \sum_{i=1}^{m-1} j_{i} \geq 1$;
%i.e., $\mathcal {T}_{m-1} = \mathcal {T}_{m-1}^{(0)}$ with $\mathcal
%{T}_{m-1}^{(0)}$ defined as
%\begin{align}
%\mathcal {T}_{m-1}^{(0)} &: = \bigg( \bigcup_{i=1}^{m-1}S_i,
%\emptyset, \emptyset, \ldots, \emptyset\bigg).
%\end{align}
%
%
%\begin{align}
%& G_{m-1} \times  \max_{\mathcal {T}_m\in
%\mathbb{T}_m}\big\{C_{\mathcal {T}_m}\big\} \nonumber
%\end{align}

%
%It's straightforward to derive
%\begin{align}
% H_{n,m}  & : = \begin{cases} 1- \frac{\prod_{i=1}^{m-1} \binom{P_n -
%i K_n}{K_n}} {\big[\binom{P_n }{K_n}\big]^{m-1}}, & \textrm{if }P_n
%\geq m K_n.
%\\ 1, & \textrm{if }P_n
%< m K_n.\end{cases} \label{eqn_Hnm},
%\end{align}

\noindent \textbf{b)} $\boldsymbol{ p_n \geq  n^{-\delta} (\ln
n)^{-1}}$

We relate $H_{n,m}$ to $H_{n,m-1}$ and assess $H_{n,m}$ iteratively.
First, with $\mathcal {T}_m^{*} = (S_1^*, S_2^*, \ldots, S_m^*)$,
event $(\mathcal {T}_m = \mathcal {T}_m^{*})$ is the intersection of
independent events: $(\mathcal {T}_{m-1} = \mathcal {T}_{m-1}^{*})$
and $(S_m = S_m^*)$. Then we have\vspace{-2pt} 
\begin{align}
  \hspace{-5pt} H_{n,m}  
& \hspace{-1pt} \vspace{-2pt} = \hspace{-1pt}\sum_{ \begin{subarray} ~\mathcal {T}_{m-1}^* \in
\mathbb{T}_{m-1} , \\  \hspace{9pt}S_m^* \in \mathbb{S}_m
\end{subarray} }\hspace{-1pt} \Big( \mathbb{P}[(\mathcal {T}_{m-1} = \mathcal
{T}_{m-1}^{*}) \mathlarger{\cap} (S_m =
S_m^*)] \times  \nonumber \\
&  \vspace{-2pt}  \quad\quad \quad\quad e^{\frac{n q_n p_n}{K_n} \sum_{1\leq i <j
\leq m-1}|S_{ij}^*|} e^{\frac{n q_n {p_n} }{K_n} \sum_{i
=1}^{m-1}|S_{i  m}^*|} \Big)  \nonumber \\ & \hspace{-1pt}  \vspace{-2pt} = \hspace{-1pt} H_{n,m-1} \cdot
\sum_{S_m^* \in \mathbb{S}_m} \mathbb{P}[ S_m  = S_m^* ]  e^{\frac{n
q_n p_n}{K_n} \sum_{i =1}^{m-1}|S_{i  m}^*|} . \label{HnmHnm1}
\end{align}
By $ \sum_{i =1}^{m-1} |S_{i  m}^*|
 \leq  m \big|S_m^*  \bcap
 \big(\bigcup_{i =1}^{m-1} S_{i }^* \big) \big|$ and (\ref{eq_pe_upper}), we have $e^{\frac{n q_n p_n }{K_n} \sum_{i =1}^{m-1}|S_{i  m}^*|} \leq e^{
\frac{2m p_n \ln n}{K_n} |S_m^* \cap
 (\bigcup_{i =1}^{m-1}S_{i }^* ) |} ,$ which is used in (\ref{HnmHnm1}) to induce\vspace{-2pt} 
\begin{align}
&  \frac{H_{n,m} }{ H_{n,m-1}}
 \leq \sum_{u=0}^{K_n} \mathbb{P}\bigg[\bigg|S_m^*  \hspace{-1pt} \bcap  \hspace{-1pt}
\bigg(\bigcup_{i =1}^{m-1}S_{i }^*\bigg)\bigg|  \hspace{-1pt} =  \hspace{-1pt} u \bigg] e^{\frac{2
u m  {p_n} \ln n}{K_n}} .\vspace{-2pt} \label{eqn_tmtm-1}
\end{align}
%
%Our goal is to further evaluate (\ref{eqn_tmtm-1_val}) based on the
%above.
Denoting $\big|\bigcup_{i=1}^{m-1}S_{i}^*\big|$ by $v$, then for $u$ satisfying
$0 \leq u \leq |S_m^*| = K_n$ and $S_m^* \bcup
\big(\bigcup_{i=1}^{m-1}S_{i}^*\big) = K_n + v - u \leq P_n $ (i.e., for $u \in [\max\{0, K_n + v - P_n\}  , K_n] $), we
obtain\vspace{-2pt} 
\begin{align}
 \hspace{-5pt} \mathbb{P}\bigg[\bigg|S_m^* \hspace{-1pt} \bcap  \hspace{-1pt}
\bigg(\bigcup_{i=1}^{m-1}S_{i}^*\bigg)\bigg|  \hspace{-1pt} =  \hspace{-1pt} u \bigg] &  =
 \binom{v}{u} \binom{P_n - v}{K_n - u} \bigg/{\binom{P_n}{K_n}},\vspace{-2pt} 
\label{probsm}
\end{align}
which together with $ K_n \leq  v \leq m K_n$ yields\vspace{-2pt} 
\begin{align}
 \textrm{L.H.S. of (\ref{probsm})} &  \leq \frac{(m K_n)^u}{u!} \hspace{-1pt} \cdot \hspace{-1pt}
 \frac{(P_n - K_n)^{K_n - u}}{(K_n - u)!} \hspace{-1pt} \cdot \hspace{-1pt} \frac{K_n !}{(P_n - K_n)^{K_n}}\vspace{-2pt} 
\nonumber
\\&   \leq \frac{1}{u!} \bigg( \frac{m {K_n}^2}{P_n - K_n}\bigg)^u. \vspace{-2pt} \label{probsm2}
\end{align}
For $u \notin [\max\{0, K_n + v - P_n\}  , K_n] $, L.H.S. of
(\ref{probsm}) equals 0. Then from (\ref{eqn_tmtm-1}) and
(\ref{probsm2}),\vspace{-2pt} 
\begin{align}
\textrm{R.H.S. of (\ref{eqn_tmtm-1})} &  \leq  \sum_{u=0}^{K_n}
\frac{1}{u!} \bigg( \frac{m {K_n}^2}{P_n - K_n} \cdot e^{\frac{2 m
{p_n} \ln n}{K_n}}\bigg)^u\vspace{-2pt}  \nonumber
\\& \quad \leq  e^{\frac{m {K_n}^2}{P_n - K_n} \cdot e^{\frac{2 m  {p_n} \ln
n}{K_n}}}. \vspace{-2pt} \label{umKnPN}
\end{align}

 By \cite[Fact 5]{ZhaoYaganGligor} and $1-x \leq e^{-x}$ for any
real $x$, it holds that\vspace{-2pt} 
\begin{align}
 s_n &  \geq 1 -  \big( 1 - K_n / P_n \big)^{K_n}  \geq 1 -   e^{- {K_n}^2/{P_n}
},\vspace{-2pt} 
 \label{eqps2}
\end{align}
For $n$ sufficiently large, from $p_n \geq n^{-\delta} (\ln n)^{-1}$
and (\ref{eq_pe_upper})
  (i.e., $q_n =p_n s_n   \leq \frac{2\ln n}{n}$), we have\vspace{-2pt} 
\begin{align}
s_n   & =  {p_n} ^{-1} {q_n} \leq {p_n} ^{-1} \cdot 2n^{-1}\ln n
\leq 2 n^{\delta-1} (\ln n)^2. \vspace{-2pt} \label{eqps0}
\end{align}
Hence, for $n$ sufficiently large, we apply (\ref{eqps2})
(\ref{eqps0}) and $P_n  >  2K_n$ (which holds from the condition $\frac{K_n}{P_n} = o (1)$) to produce\vspace{-2pt} 
\begin{align}
  & {{K_n}^2}/({P_n-K_n}) <  {2 {K_n}^2}/{P_n} \leq  -2\ln (1 -
s_n)\vspace{-2pt}   \nonumber
\\& \quad \leq  -2\ln (1 - 2 n^{\delta-1} (\ln n)^2) \leq 2\sqrt{2 } n^{\frac{\delta-1}{2}} \ln n, \vspace{-2pt} \label{eqn_knpn2}
\end{align}
where the last step uses $-\ln (1-y) \leq \sqrt{y}$ for $0<y<1$. From (\ref{eq_pe_lnnn}) and condition $P_n = \Omega (n)$, we obtain from \cite[Lemma 7]{ZhaoYaganGligor} that $K_n =  \omega\big( \sqrt{\ln n} \big) =  \omega(1)  $.
Then for an arbitrary constant $c > 2$, it holds that $\frac{K_n}{p_n} \geq K_n \geq \frac{4c \cdot
m}{(c-2)(1-\delta)} $ holds for
all $n$ sufficiently large. Hence,\vspace{-2pt}
\begin{align}
e^{\frac{2 m p_n \ln n}{K_n}} & \leq e^{  \frac{(c-2)(1-\delta)}{2c}
\ln n} = n^{\frac{(c-2)(1-\delta)}{2c}} .\label{ja1}\vspace{-2pt}
\end{align}
The use of (\ref{umKnPN}) (\ref{eqn_knpn2}) and (\ref{ja1}) in
(\ref{eqn_tmtm-1}) yields\vspace{-2pt}
\begin{align}
 & H_{n,m} / H_{n,m-1} \leq \textrm{R.H.S. of (\ref{eqn_tmtm-1})}
\vspace{-2pt} \nonumber \\ & \quad \leq  e^{ 2\sqrt{2 } m n^{\frac{\delta-1}{2}} \cdot
n^{\frac{(c-2)(1-\delta)}{2c}} \cdot  \ln n } \leq \Big(e^{3
n^{\frac{\delta-1}{c}} \ln n} \Big)^m.\vspace{-2pt} \label{gnmgnm-1}
\end{align}

To derive $H_{n,m}$ iteratively based on (\ref{gnmgnm-1}), we
compute $H_{n,2}$ below. Setting $m=2$ in L.H.S. of
(\ref{eqn_sumTmst}) and considering the independence between  
$(S_1  = S_1^*)$ and $(S_2  = S_2^*)$, we gain\vspace{-2pt}
\begin{align}
\hspace{-2pt}  H_{n,2} &  \hspace{-2pt} =  \hspace{-4pt} \sum_{S_1^*
\in \mathbb{S}_m} \hspace{-3pt}  \mathbb{P}[ S_1  = S_1^* ]
 \hspace{-4pt}   \sum_{S_2^* \in \mathbb{S}_m} \hspace{-2pt}  \mathbb{P}[ S_2  =
S_2^* ] e^{\frac{n q_n p_n}{K_n} |S_1^* \cap S_2^*|}.\vspace{-2pt}
\label{eqn_gn2}
\end{align}
Clearly, $\sum_{S_2^* \in \mathbb{S}_m} \hspace{-3pt} \mathbb{P}[
S_2 \hspace{-1pt} = \hspace{-1pt} S_2^* ] e^{\frac{n q_n p_n}{K_n}
|S_1^* \cap S_2^*|} $ equals R.H.S. of (\ref{eqn_tmtm-1}) with $m =
2$. Then from (\ref{gnmgnm-1}) and (\ref{eqn_gn2}),\vspace{-2pt}
\begin{align}
H_{n,2}  &  \leq \sum_{S_1^* \in \mathbb{S}_m} \mathbb{P}[ S_1  =
S_1^* ]  e^{6 n^{\frac{\delta-1}{c}} \ln n} =   e^{6
n^{\frac{\delta-1}{c}} \ln n}. \vspace{-2pt}\label{hn2}
\end{align}

Therefore, it holds via (\ref{gnmgnm-1}) and (\ref{hn2}) that\vspace{-2pt}
\begin{align}
H_{n,m}  \hspace{-2pt}  \leq \hspace{-2pt} \big(\hspace{-1.5pt}e^{3 n^{\frac{\delta-1}{c}} \ln n}
\hspace{-1.5pt}\big)^{\hspace{-1.5pt}m+(m-1) + \ldots + 3}  \hspace{-1.5pt} e^{6  n^{\frac{\delta-1}{c}} \hspace{-1pt}
\ln n} \hspace{-2pt} \leq \hspace{-2pt} e^{3m^2 n^{\frac{\delta-1}{c}} \ln n}\vspace{-2pt}
\nonumber.
\end{align}

Finally, from cases a) and b), for $n$ sufficiently large, $ H_{n,m} $ is at most $ \max\big\{e^{ m^2 n^{-\delta}} , \hspace{2pt} e^{3m^2 n^{\frac{\delta-1}{c}} \ln n}\big\}$.
%\begin{align}
% H_{n,m} & \leq  \max\left\{e^{ m^2 n^{-\delta}} , e^{3m^2 n^{\frac{\delta-1}{c}} \ln n}\right\} .
%\nonumber
%\end{align}
Then (\ref{eqn_sumTmst}) follows.\vspace{-2pt}

\section{The Proof of Proposition 2} \label{sec:PROP_SND}\vspace{-2pt}

%
%\begin{align}
%&  \mathbb{P} \big[ \big( \mathcal {L}_m = \mathcal {L}_m^{*} \big)
%\cap \big( \mathcal {M}_m = \mathcal{M}_m^{(0)}
% \big) \big]  \nonumber  \\
%& \quad =  \mathbb{P} \big[ \mathcal {L}_m = \mathcal {L}_m^{*}
% \big] \mathbb{P} \big[ \big( \mathcal {M}_m = \mathcal{M}_m^{(0)}
% \big)  \boldsymbol{\mid} \big( \mathcal {L}_m = \mathcal
%{L}_m^{*} \big) \big] \nonumber  \\
%& \quad = \mathbb{P} \big[ \mathcal {C}_m = \mathcal {C}_m^{*} \big]
%\mathbb{P} \big[ \mathcal {T}_m = \mathcal {T}_m^{*}
% \big]  \mathbb{P} \big[ \big( \mathcal {M}_m = \mathcal{M}_m^{(0)} \big)  \boldsymbol{\mid} \big( \mathcal {T}_m =
%\mathcal {T}_m^{*}\big) \big]\nonumber
%\end{align}
%
%For any $\mathcal {L}_m^{*} \in \mathbb{L}_m^{(0)}$,
%\begin{align}
%\mathbb{M}_m (\mathcal {L}_m^{*}) = \mathbb{M}_m^{(0)} \nonumber
%\end{align}
%
%\begin{align}
%&  \mathbb{P} \big[ \big( \mathcal {M}_m = \mathcal {M}_m^{*}
%(\mathcal {L}_m^{*}) \big)  \boldsymbol{\mid} \big( \mathcal {L}_m =
%\mathcal
%{L}_m^{*} \big) \big] \nonumber  \\
%& \quad =  \mathbb{P} \big[ \big( \mathcal {M}_m = \mathcal
%{M}_m^{*} (\mathcal {L}_m^{*}) \big)  \boldsymbol{\mid} \big(
%\mathcal {T}_m = \mathcal {T}_m^{*}\big) \big]\nonumber
%\end{align}
%
%\begin{align}
%\mathcal {L}_m^{*} & = (\mathcal {C}_m^{*}, \mathcal {T}_m^{*})
%\nonumber
%\end{align}

We define \vspace{-2pt}$\mathcal{C}_m^{(0)}$ \vspace{-1pt}and $\mathbb{T}_m^{(0)}$ by $\mathcal {C}_m^{(0)}  = ( \underbrace{0, 0, \ldots,
0}_{m(m-1)/2 \textrm{ number of ``}0\textrm{''}} )$
%\begin{align}
%\mathcal {C}_m^{(0)}  = ( \underbrace{0, 0, \ldots,
%0}_{\binom{m}{2} \textrm{ number of ``}0\textrm{''}} ), \nonumber
%\end{align}
and $\mathbb{T}_m^{(0)} \hspace{-2pt}
  = \hspace{-2pt}
 \{\mathcal {T}_m \hspace{-2pt}
 \boldsymbol{\mid} \hspace{-2pt}
 S_i \cap
S_j = \emptyset, ~ \forall  1 \leq i < j \leq
m.\}$.
%\begin{align}
%\mathbb{T}_m^{(0)} \hspace{-2pt}
%  = \hspace{-2pt}
% \{\mathcal {T}_m \hspace{-2pt}
% \boldsymbol{\mid} \hspace{-2pt}
% S_i \cap
%S_j = \emptyset, ~ \forall i, j\textrm{ with }1 \leq i < j \leq
%m.\}. \nonumber
%\end{align}
Clearly, $\big(\mathcal {C}_m \hspace{-2pt}=\hspace{-2pt} \mathcal{C}_m^{(0)}\big)$ or
$\big(\mathcal {T}_m \hspace{-2pt}\in\hspace{-2pt} \mathbb{T}_m^{(0)}\big)$ each implies
$\big( \mathcal {L}_m \hspace{-2pt}\in\hspace{-2pt} \mathbb{L}_m^{(0)} \big)$. Also,
$\big(\mathcal {C}_m \hspace{-1pt} = \hspace{-1pt} \mathcal{C}_m^{(0)}\big)$ and $\big(\mathcal
{M}_m \hspace{-1pt} = \hspace{-1pt} \mathcal{M}_m^{(0)}\big)$ are independent of each other.
Thus, with $\mathcal {P}_2 \hspace{-1pt}=\hspace{-1pt} \mathbb{P} \big[ \big( \mathcal
{L}_m \in \mathbb{L}_m^{(0)} \big) \cap \big( \mathcal {M}_m =
\mathcal{M}_m^{(0)} \big) \big]$, \vspace{-2pt}we derive
\begin{align}
& \mathcal {P}_2 \geq \mathbb{P} \big[ \mathcal {C}_m =
\mathcal{C}_m^{(0)}\big] \mathbb{P} \big[ \mathcal {M}_m =
\mathcal{M}_m^{(0)} \big], \vspace{-2pt}\label{prcm}
\end{align}
and\vspace{-2pt}
\begin{align}
& \mathcal {P}_2 \geq \mathbb{P} \big[ \mathcal {T}_m \hspace{-1pt}
\in  \hspace{-1pt} \mathbb{T}_m^{(0)}
 \big] \mathbb{P} \big[ \big( \mathcal {M}_m  \hspace{-1pt} =  \hspace{-1pt} \mathcal{M}_m^{(0)} \big) \hspace{-2pt}
  \boldsymbol{\mid} \hspace{-2pt}
 \big( \mathcal {T}_m  \hspace{-1pt} \in
\mathbb{T}_m^{(0)}  \hspace{-1pt} \big)\big].\vspace{-2pt} \label{eqn_tmtmst}
\end{align}

Given that event $\big(\mathcal {C}_m = \mathcal{C}_m^{(0)}\big) $ is $
\overline{\bigcup_{ 1 \leq i < j \leq m} {C_{ij}}} $ and
event $\big(\mathcal {T}_m \in \mathbb{T}_m^{(0)}\big)  $ is $ \overline{\bigcup_{ 1
\leq i < j \leq m} {\Gamma_{ij}} }$, \vspace{-2pt}using the union bound, we
get
\begin{align}
 & \mathbb{P} \big[ \mathcal {C}_m = \mathcal{C}_m^{(0)} \big]\geq 1 -
\sum_{ 1 \leq i < j \leq m}\mathbb{P}[ C_{ij} ] \geq 1- m^2 p_n / 2,\vspace{-2pt}
\label{prcmpn}
\end{align}
and
\begin{align}
& \mathbb{P}\big[\mathcal {T}_m  \in \mathbb{T}_m^{(0)}\big] \geq 1
-\vspace{-2pt} \sum_{ 1 \leq i < j \leq m}\mathbb{P}[ \Gamma_{ij} ] \geq  1 - m^2
s_n / 2.\label{mthbbP}
 % \label{eqn_tm_tmstar}
\end{align}

Denoting $(h!)^{-m} (n q_n)^{hm} e^{-m n q_n}$ by $\Lambda$,
we will prove\vspace{-2pt}
\begin{align}
\mathbb{P}  \big[ \mathcal {M}_m =\vspace{-2pt} \mathcal{M}_m^{(0)} \big] & \sim
\Lambda \label{eqn_prMm},
\end{align}
and\vspace{-2pt}
\begin{align}
&  \mathbb{P} \big[ \big( \mathcal {M}_m = \mathcal{M}_m^{(0)} \big)
\boldsymbol{\mid} \big( \mathcal {T}_m \in \mathbb{T}_m^{(0)}
\big)\big]   \vspace{-2pt}\geq  \Lambda \cdot [1-o(1)]
.\label{prob_MmMm_sim}
\end{align}
Substituting (\ref{prcmpn}) and (\ref{eqn_prMm}) into (\ref{prcm}),
and applying (\ref{mthbbP}) and (\ref{prob_MmMm_sim}) to
(\ref{eqn_tmtmst}), we get (i) $\mathcal {P}_2 / \Lambda   \geq (  1 - \min\{ s_n, p_n \}  \cdot m^2  / 2)  [1-o(1)]
  . $\\
%\begin{align}
%&  \mathcal {P}_2 / \Lambda   \geq (  1 - \min\{ s_n, p_n \}  \cdot m^2  / 2)\cdot [1-o(1)]
%  .\vspace{-2pt} \label{pro2_pt1}
%\end{align}
From (\ref{eqn_prMm}), we get (ii) $ \mathcal {P}_2  \hspace{-2pt} \leq \hspace{-2pt}\mathbb{P} \big[ \mathcal {M}_m\hspace{-2pt}
\in\hspace{-2pt}
\mathbb{M}_m^{(0)} \big] \hspace{-2pt} \leq \hspace{-2pt} \Lambda   [1\hspace{-2pt}+\hspace{-2pt}o(1)].$
\vspace{-2pt}
%\begin{align}
% \mathcal {P}_2   \leq \mathbb{P} \big[ \mathcal {M}_m
%\in
%\mathbb{M}_m^{(0)} \big]  \leq  \Lambda \cdot [1+o(1)]. \vspace{-2pt}
%\label{pro2_pr}
%\end{align}
Combining (i) and (ii) above %(\ref{pro2_pt1}) and (\ref{pro2_pr}), 
and using\vspace{1pt} $\min\{
s_n, p_n \} \leq \sqrt{s_n p_n}  = \sqrt{q_n} = o(1)$ which holds from $q_n = s_n p_n$ and
(\ref{eq_pe_lnnn}), Proposition 2 follows. Below we establish (\ref{eqn_prMm}) and (\ref{prob_MmMm_sim}).\vspace{-2pt}

%
%\begin{align}
%& \mathbb{P} \big[ \big( \mathcal {L}_m \in \mathbb{L}_m^{(0)} \big)
%\cap \big( \mathcal {M}_m = \mathcal{M}_m^{(0)} \big) \big]
% \nonumber  \\
%& \quad \geq \bigg[ 1- \frac{m(m-1)}{2} s_n \bigg] \cdot \mathbb{P}
%\big[ \big( \mathcal {M}_m = \mathcal{M}_m^{(0)} \big)
%\boldsymbol{\mid} \big( \mathcal {T}_m \in \mathbb{T}_m^{(0)}
%\big)\big]
%\end{align}

\subsection{Establishing (\ref{eqn_prMm})}\vspace{-2pt}

We write $\mathbb{P} \big[ \mathcal {M}_m = \mathcal{M}_m^{(0)}
\big]$ as \vspace{-2pt}
\begin{align}
& \sum_{\mathcal {T}_m^{*} \in \mathbb{T}_m} \hspace{-3pt} \Big\{
\mathbb{P} \big[ \mathcal {T}_m \hspace{-2pt} = \hspace{-2pt}
\mathcal {T}_m^{*} \big] \mathbb{P} \big[ \big( \mathcal {M}_m
\hspace{-2pt} = \hspace{-2pt} \mathcal{M}_m^{(0)} \big)
\boldsymbol{\mid} \big( \mathcal {T}_m \hspace{-2pt} = \hspace{-2pt}
\mathcal {T}_m^{*} \big)\big] \Big\},\nonumber \vspace{-2pt}
\end{align}
where $\mathbb{P} \big[ \big( \mathcal {M}_m = \mathcal{M}_m^{(0)}
 \big) \boldsymbol{\mid} \big( \mathcal {T}_m
= \mathcal {T}_m^{*} \big)\big] $ equals
\begin{align}\vspace{-2pt}
&  f\big(n-m , \mathcal{M}_m^{(0)}\big) \mathbb{P} [w \in M_{0^m}
\boldsymbol{\mid}%(\mathcal {M}_m = \mathcal
%{M}_m^{*}) \cap
\mathcal {T}_m = \mathcal {T}_m^{*} ]^{n-m-hm} \nonumber  \\
& \quad \quad \times \prod_{i=1}^{m} \mathbb{P}[w \in
M_{0^{i-1}, 1, 0^{m-i}} \boldsymbol{\mid}%(\mathcal {M}_m = \mathcal
%{M}_m^{*}) \cap
\mathcal {T}_m = \mathcal {T}_m^{*} ]^{h}, \vspace{-2pt} \nonumber
\end{align}
where $f\big(n-m , \mathcal{M}_m^{(0)}\big)$ is the number of ways assigning the $(n-m)$ nodes from $ \overline{\mathcal {V}_m} $ to $M_{j_1 j_2 \ldots j_m}$ such that       $|M_{j_1 j_2 \ldots j_m}|$ is given by $\mathcal{M}_m^{(0)}$ (see (\ref{mmsizesat})). Hence, it holds from (\ref{eqn_fnexpr}) that
\vspace{-2pt}
\begin{align}
& \hspace{-2pt}  f\big(n\hspace{-1pt}-\hspace{-1pt}m ,
\mathcal{M}_m^{(0)} \hspace{-1pt}\big) \hspace{-2pt} \vspace{-2pt}= \hspace{-2pt}
\frac{ (n \hspace{-1pt}- \hspace{-1pt}m) ! }{(n
\hspace{-1pt}-\hspace{-1pt} m\hspace{-1pt} -\hspace{-1pt}
hm)!(h!)^m} \hspace{-3pt} \sim \hspace{-2pt} (h!)^{-m}n^{hm}.
\label{eqn_f00}
\end{align}
 We will establish \vspace{-2pt}
\begin{align}
 &\hspace{-4pt}
\vspace{-2pt}  \sum_{ \mathcal {T}_m^{*} \in \mathbb{T}_m }
   \hspace{-3pt}\Big\{  
\mathbb{P}[\mathcal {T}_m \hspace{-1pt}
 = \hspace{-1pt}
 \mathcal {T}_m^{*}] \hspace{1pt}
\prod_{i=1}^{m} \{ \mathbb{P}\big[w \hspace{-2pt}
 \in \hspace{-2pt}
 M_{0^{i-1}, 1,
0^{m-i}}   \hspace{-2pt}\boldsymbol{\mid} \hspace{-2pt}
\mathcal {T}_m \hspace{-2pt}
 = \hspace{-2pt}
 \mathcal {T}_m^{*} \big]^h  \} \Big\} 
\nonumber  \\
& \quad  \geq {q_n}^{hm} \cdot [1-o(1)]  \vspace{-2pt} .\label{eq_evalprob_exp_2}
\end{align}
We use (\ref{eqn_f00}) and (\ref{eq_evalprob_exp_2}) as well as
(\ref{eq_evalprob_3}) (viz., Lemma \ref{lem_evalprob} in the Appendix) in evaluating
$\mathbb{P} \big[ \mathcal {M}_m = \mathcal{M}_m^{(0)} \big]$ above.
Then \vspace{-2pt}
\begin{align}
&  \mathbb{P} \big[ \mathcal {M}_m = \mathcal{M}_m^{(0)} \big]\vspace{-2pt}
\nonumber  \\
& \geq  (h!)^{-m}n^{hm} \cdot [1-o(1)] \cdot (1-m q_n)^{n} \times \vspace{-2pt}
 \nonumber  \\
& \quad \sum_{\mathcal {T}_m^{*} \in \mathbb{T}_m}
\hspace{-2pt} \mathbb{P}[\mathcal {T}_m \hspace{-1pt}
 = \hspace{-1pt}
 \mathcal {T}_m^{*}] \hspace{1pt}
\prod_{i=1}^{m}\big\{ \mathbb{P}[w  \hspace{-2pt}  \in \hspace{-2pt}
 M_{0^{i-1},
1, 0^{m-i}} \hspace{-2pt}
 \boldsymbol{\mid} \hspace{-2pt}
 \mathcal {T}_m \hspace{-2pt}
 =  \hspace{-2pt} \vspace{-2pt}
 \mathcal {T}_m^{*} ]^{h} \big\} \nonumber \\
& \geq  (h!)^{-m} (n q_n)^{hm} e^{-m n q_n} \cdot [1-o(1)] .\vspace{-2pt}
\vspace{-2pt} \label{eqn_prMm_pt2}
\end{align}

Substituting (\ref{EQ}) (\ref{eqn_f00}) above
 and (\ref{eq_evalprob_1}) in Lemma \ref{lem_evalprob} into the computation of $\mathbb{P} \big[ \mathcal {M}_m = \mathcal{M}_m^{(0)}
\big]$ yields \vspace{-2pt}
\begin{align}
& \mathbb{P} \big[ \mathcal {M}_m = \mathcal{M}_m^{(0)} \big]\vspace{-2pt}
\nonumber  \\
& \leq  (h!)^{-m}n^{hm}\vspace{-2pt} {q_n}^{hm} \times [1+o(1)] \times \nonumber  \\
&\quad \sum_{\mathcal {T}_m^{*} \in \mathbb{T}_m} \mathbb{P}[\mathcal
{T}_m = \mathcal {T}_m^{*}] \vspace{-2pt}\hspace{1pt}\mathbb{P} [w \in M_{0^m}
\boldsymbol{\mid} \mathcal {T}_m = \mathcal {T}_m^{*}
]^{n-m-hm}  \nonumber \\
& \sim (h!)^{-m} (n q_n)^{hm} e^{-m n q_n} .\vspace{-2pt} \label{eqn_prMm_pt1}
\end{align}

Then (\ref{eqn_prMm}) follows from (\ref{eqn_prMm_pt2}) and
(\ref{eqn_prMm_pt1}). Namely, (\ref{eqn_prMm}) holds upon the
establishment of (\ref{eq_evalprob_exp_2}).
\iffalse
First, from (\ref{eq_evalprob_2}) in Lemma \ref{lem_evalprob}, with $\mathcal
{T}_m^{*} = (S_1^{*} , S_2^{*}  , \ldots, S_m^{*} ) $ and
$S_{ij}^{*} = S_{i}^{*} \cap S_{j}^{*}$, we get
\begin{align}
& \prod_{i=1}^{m} \mathbb{P}\big[w \in M_{0^{i-1}, 1, 0^{m-i}} 
\boldsymbol{\mid} \mathcal {T}_m = \mathcal {T}_m^{*} \big]^h
\nonumber  \\
& \geq { {q_n}^{hm} \hspace{-1pt} \prod_{i=1}^{m} \bigg[ 1
\hspace{-1pt} - \hspace{-2pt} \bigg( \hspace{-2pt} 2m{q_n}
\hspace{-1pt}
 + \hspace{-1pt} \frac{p_n}{K_n} \hspace{-2pt}
\sum_{j\in\{1,2,\ldots,m\}\setminus\{i\}}
 |S_{ij}^{*}|\bigg)\hspace{-2pt} \bigg] }^h\nonumber  \\
&  \geq {q_n}^{hm} \bigg( 1 - 2 h m^2 {q_n}  - \frac{2 hm p_n}{K_n}
\sum _{1\leq i <j \leq m}  |S_{ij}^{*}|\bigg), 
\end{align}
where the last step uses the following inequality easily proved by mathematical induction: $\prod_{\ell=1}^{r} (1-x_{\ell}) \geq 1-\sum_{\ell=1}^r x_{\ell} $ for any positive integer $r$ and any positive $ x_{\ell}$ for $\ell=1,2,\ldots, r$ (we set $\ell = mh$, with the $mh$ number of $x_l$ having $m$ groups, where the group $i$ for $i=1,2,\ldots, m$ has $m$ members all being $2m{q_n}
\hspace{-1pt}
 + \hspace{-1pt} \frac{p_n}{K_n} \hspace{-2pt}
\sum_{j\in\{1,2,\ldots,m\}\setminus\{i\}}
 |S_{ij}^{*}|$.)
 \fi
 From (\ref{hm2qnpnleq})  in Lemma \ref{lem_evalprob} and $q_n = o(1)$ by (\ref{eq_pe_lnnn}), we obtain
(\ref{eq_evalprob_exp_2}) once proving\vspace{-2pt}
\begin{align}
\frac{ p_n}{K_n} \hspace{-2pt} \sum_{ \mathcal {T}_m^{*} \in
\mathbb{T}_m } \hspace{-2pt} \hspace{-2pt} \bigg(
\mathbb{P}[\mathcal {T}_m = \mathcal {T}_m^{*}] \hspace{-2pt} \sum
_{1\leq i <j \leq m} \hspace{-2pt}|S_{ij}^{*}| \hspace{-1pt} \bigg)
& = o(1). \hspace{-2pt}\vspace{-2pt} \label{prfhpnkn}
\end{align}
If $\mathcal {T}_m^{*}   \in \mathbb{T}_m^{(0)}$, then $| S_{ij}^{*}
|  =   0  $. Then from (\ref{mthbbP}), we get 
(\ref{prfhpnkn}) by\vspace{-2pt}
\begin{align}
 \textrm{L.H.S. of (\ref{prfhpnkn})}
  & \leq p_n \hspace{-1pt} \cdot \hspace{-1pt} m(m-1) /2 \hspace{-1pt}  \cdot \hspace{-1pt} \mathbb{P}[\mathcal
{T}_m^{*} \in \mathbb{T}_m \setminus
\mathbb{T}_m^{(0)}] \vspace{-2pt}\nonumber  \\
&   \leq p_n \hspace{-1pt} \cdot \hspace{-1pt} m^2 /2  \hspace{-1pt} \cdot \hspace{-1pt} m^2 s_n / 2 \leq m^4
n^{-1}\ln n / 2 \vspace{-2pt} = o(1).
 \nonumber
\end{align}

\subsection{Establishing (\ref{prob_MmMm_sim})}\vspace{-2pt}

%By the definition of $\mathbb{T}_m^{(0)}$,

Let $\Delta$ denote $\mathbb{P} \big[ \big( \mathcal {M}_m \hspace{-2pt} = \hspace{-2pt}
\mathcal{M}_m^{(0)} \big) \boldsymbol{\mid} \big( \mathcal {T}_m
\hspace{-2pt} \in \hspace{-2pt} \mathbb{T}_m^{(0)} \big)\big]$. Clearly,
$\Delta$ is equivalent to $\mathbb{P} \big[ \big(
\mathcal {M}_m = \mathcal{M}_m^{(0)}
 \big) \boldsymbol{\mid} \big( \mathcal {T}_m
= \mathcal {T}_m^{*} \big)\big]$ for any $\mathcal {T}_m^{*} \in
\mathbb{T}_m^{(0)}$, so it follows that \vspace{-2pt}
\begin{align}
\hspace{-6pt} \Delta & \hspace{-2pt} = \hspace{-2pt}
f\big(n\hspace{-1pt}-\hspace{-1pt}m , \mathcal{M}_m^{(0)}\big)
\mathbb{P} [w \in M_{0^m}
\hspace{-2pt}\boldsymbol{\mid}\hspace{-2pt}%(\mathcal {M}_m = \mathcal
%{M}_m^{*}) \cap
\mathcal {T}_m \hspace{-1pt}=\hspace{-1pt} \mathcal {T}_m^{*} ]^{n-m-hm}\vspace{-2pt}\nonumber  \\
& ~~\times \prod_{i=1}^{m}\big\{ \mathbb{P}[w \in
M_{0^{i-1}, 1, 0^{m-i}} \boldsymbol{\mid}%(\mathcal {M}_m = \mathcal
%{M}_m^{*}) \cap
\mathcal {T}_m \vspace{-2pt}= \mathcal {T}_m^{*} ]^{h} \big\},\label{eqn_probMm}
\end{align}
with $f\big(n-m , \mathcal{M}_m^{(0)}\big)$ given by
(\ref{eqn_f00}). For $\mathcal {T}_m^{*} \in \mathbb{T}_m^{(0)}$,
from $|S_{ij}^{*}| = 0$ and (\ref{hm2qnpnleq}) in Lemma \ref{lem_evalprob}, we
derive\vspace{-2pt}
\begin{align}
\prod_{i=1}^{m}  \hspace{-1pt}\Big\{ \mathbb{P}\big[w \hspace{-1pt} \in \hspace{-1pt} M_{0^{i-1}, 1, 0^{m-i}} 
\hspace{-1pt} \boldsymbol{\mid} \hspace{-1pt} \mathcal {T}_m \hspace{-1pt} =\hspace{-1pt}  \mathcal {T}_m^{*} \big] \Big\}^h &\hspace{-1pt} \vspace{-2pt} \geq\hspace{-1pt}  {q_n}^{hm} ( 1\hspace{-1pt}  -\hspace{-1pt}  2 h m^2 {q_n}  )
\label{eqn_wM0} . \vspace{-2pt}
\end{align}
Substituting (\ref{eqn_f00}) (\ref{eqn_wM0}) above and
(\ref{eq_evalprob_3}) in Lemma \ref{lem_evalprob} into (\ref{eqn_probMm}), we
conclude that $\Delta$ is at least \vspace{-2pt}
\begin{align}
&    (h!)^{-m}n^{hm} \cdot [1-o(1)]
\vspace{-2pt}  \nonumber  \\
 &  \times  {q_n}^{hm} ( 1 - 2 h m^2 {q_n}  )
\hspace{-1pt}   \cdot \hspace{-1pt}  (1 \hspace{-1pt} - \hspace{-1pt} m q_n )^{n-m-hm} \hspace{-1pt}\vspace{-2pt} =\hspace{-1pt}   \Lambda  \cdot [1\hspace{-1pt} -\hspace{-1pt} o(1)]. \nonumber
\end{align}

\section{Numerical Experiments} \label{sec:expe}

To confirm our analytical results, we now provide numerical
experiments in the non-asymptotic regime. %As we will see from the simulation, the experimental
%observations are in agreement with
% our theoretical findings.

\begin{figure}[!t]
  \centering
  \psfrag{Z}{\centerhack{$\mathsmaller{\mathbb{P}[\mathbb{G}(n,K,P,p)\textrm{ is $2$-connected}.]}$}}
 \includegraphics[height=0.21\textwidth]{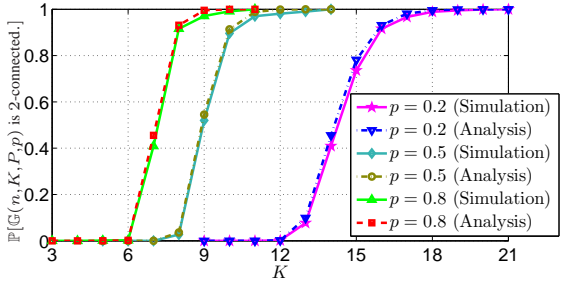}  \vspace{-5pt}
 \caption{A plot generated from the simulation and
         the analysis for the
probability that $\mathbb{G}(n,K,P,p)$
        is $2$-connected versus
$K$ with $n=2,000$, $P=10,000$ and $p =0.2, 0.5,0.8$. \vspace{-15pt} }
\label{figure:k_2_on_off}
\end{figure}

In Figure \ref{figure:k_2_on_off}, we depict the probability that
graph $\mathbb{G}(n,K,P,p)$ is $2$-connected from both the
simulation and the analysis, as elaborated below. In all set of
experiments,
%contributing to Figure \ref{figure:k_2_on_off},
 we fix the number of
nodes at $n=2,000$ and the key pool size at $P=10,000$. For the
probability $p$ of a communication channel being {\em on}, we
consider $p =0.2, 0.5,0.8$, while varying the parameter $K$ from $3$
to $21$. % ($K$ is the number of keys on each sensor).
 For each pair
$(K, p)$, we generate $1,000$ independent samples of
$\mathbb{G}(n,K,P,p)$ and count the number of times that the
obtained graphs are $2$-connected. Then the counts divided by $1,000$ become the empirical
probabilities.

%
%\begin{figure}[!t]
%  \centering
% \includegraphics[width=0.41\textwidth]{figures/f3.eps}
% \caption{\sl A plot for the
%probability that $\mathbb{G}(n,K,P,p)$ has a minimum node degree no
%less than $k$ as a function of $p$ for
%         $k=2$, $k=3$ and $k=8$
%         with $n=2,000$, $K=37$ and $P=2,000$.}  \label{figure:k_2_on_off2}
%\end{figure}

The curves in Figure \ref{figure:k_2_on_off} corresponding to the
analysis are determined as follows. We use the asymptotical result
to approximate the probability of $2$-connectivity in
$\mathbb{G}(n,K,P,p)$; specifically, given $n,K,  P, p$ and
$k=2$, we   determine $\alpha$ by considering $p \cdot
 \big[1-  {\binom{P - K}{K} } \big/ {\binom{P}{K}}\big] =
\frac{\ln n + {(k-1)} \ln \ln n + {\alpha }}{n}$, a condition
stemming from (\ref{thm_eq_pe}) and the computation of ${q_n}$ in
Section \ref{sec:SystemModel}, and then use $e^{-
\frac{e^{-\alpha}}{(k-1)!}}$ as the analytical reference of
$\mathbb{P}[\mathbb{G}(n,K,P,p)\textrm{ is 2-connected}]$ for a
comparison with the empirical probabilities. Figure
\ref{figure:k_2_on_off} indicates
  that the experimental results are
 in agreement with our analysis.

\section{Related Work} \label{related} \vspace{-3pt}

\textbf{Random key graphs.} For a random key graph $G(n,K_n,P_n)$
(viz., Section \ref{sec:SystemModel}) which models the topology
induced by the EG scheme, Rybarczyk \cite{ryb3}
derives the asymptotically exact probability of connectivity,
covering a weaker form of the result -- a zero-one law which is also
 obtained in \cite{r1,yagan}. %As demonstrated in \cite{ryb3}, in
%$G(n,K_n,P_n)$ with $K_n \geq 2$, $\frac{{K_n}^2}{P_n}=\frac{\ln n +
%{\alpha_n}}{n}$ and $\lim_{n \to \infty} \alpha_n = \alpha ^* $, the
%probability of connectivity asymptotically converges to $e^{-
%e^{-\alpha ^*}}$.
 Rybarczyk \cite{zz} further establishes a zero-one
law for $k$-connectivity, and we \cite{ZhaoCDC} obtain the
asymptotically exact probability of $k$-connectivity. Under $P_n =
\Theta(n ^{c})$ for some constant $c
>1$ and $\frac{{K_n}^2}{P_n}=\frac{\ln n + (k-1) \ln \ln n +
{\alpha_n}}{n}$, Rybarczyk's result \cite{zz} is that the
probability of $k$-connectivity in graph $G(n,K_n,P_n)$ is
asymptotically converges to $1$ (resp. $0$) if $ \lim_{n \to \infty}
\alpha_n $ equals $\infty$ (resp., $-\infty$), while we
\cite{ZhaoCDC} prove that such probability asymptotically approaches
to $e^{- \frac{e^{-\alpha ^*}}{(k-1)!}}$ if $\lim_{n \to \infty}
\alpha_n = \alpha ^* \in (-\infty, \infty)$. \iffalse As mentioned, our
Theorem
 \ref{THM1} and Corollary \ref{COR1} also imply the results on the asymptotically exact
probability of $k$-connectivity in graph $G(n,K_n,P_n)$, by setting
$p_n$ (the failure probability of a channel) as $1$. \fi

 %As noted after Theorem
%\ref{THM1}, our results also imply zero-one laws for the two
%properties in graph $G(n,K_n,P_n)$.
%
%
%\emph{almost surely}\footnote{An event relating with $n$ happens
%\textit{almost surely} if its probability converges to 1 as $n\to
%\infty$.} (a.s.)

%{\bf How about S. Blackburn and Gerke [1]  We need one line}

\textbf{Erd\H{o}s--R\'{e}nyi graphs.} For an Erd\H{o}s--R\'{e}nyi
graph $G(n,p_n)$ where any two nodes have an edge in between
independently with probability $p_n$, Erd\H{o}s and R\'{e}nyi
consider connectivity in  
\cite{citeulike:4012374} and $k$-connectivity in \cite{erdos61conn},
where the latter result is that if $p_n\hspace{-2pt} =\hspace{-2pt}
\frac{\ln n + {(k-1)} \ln \ln n + {\alpha_n}}{n}$ and $\lim_{n \to
\infty} \alpha_n\hspace{-2pt} =\hspace{-2pt} \alpha ^* \in [-\infty,
\infty] $, graph $G(n,p_n)$ is $k$-connected with a probability
asymptotically tending to $e^{- \frac{e^{-\alpha ^*}}{(k-1)!}}$.

\textbf{Random key graphs $\mathlarger{\cap}$ Erd\H{o}s--R\'{e}nyi
graphs.} As given in Section \ref{sec:SystemModel}, our studied
graph $\mathbb{G}$ is the intersection of a random key graph
$G(n,K_n,P_n)$ and an Erd\H{o}s--R\'{e}nyi graph $G(n,p_n)$. For
graph $\mathbb{G}$, Ya\u{g}an \cite{yagan_onoff} establishes a
zero-one law for connectivity, and we \cite{ZhaoYaganGligor,ISIT}
extend Ya\u{g}an's result to $k$-connectivity and show that with
$P_n = \Omega (n)$, $\frac{K_n}{P_n} = o (1)$ and $q_n$ set as
$\frac{\ln n + (k-1) \ln \ln n + {\alpha_n}}{n}$, graph $\mathbb{G}$
is (resp., is not) $k$-connected with high probability if
$\lim_{n\to \infty}{\alpha_n} = \infty$ (resp., $\lim_{n\to
\infty}{\alpha_n} = -\infty$). Compared               with this result in
\cite{ZhaoYaganGligor,ISIT}, our result on the asymptotically exact
probability of $k$-connectivity is stronger and more challenging to
derive.

\textbf{Random key graphs $\mathlarger{\cap}$ random geometric
graphs.} Connectivity properties have also been studied in secure
sensor networks employing the EG scheme under the disk model, where
any two nodes need to be within a certain distance $r_n$ to have a
link in between. When nodes are assumed to be uniformly and
independently deployed in some region $\mathcal {A}$, the topology of
such a network is represented by the intersection of a random key
graph $G(n,K_n,P_n)$ and a random geometric graph, where a random
geometric
 graph denoted by
$G(n,r_n,\mathcal {A})$ is defined on $n$ nodes independently and
uniformly
  distributed in $\mathcal{A}$ such that an edge
  exists between two nodes if and only if their distance is at most
  $r_n$. Krzywdzi\'{n}ski and Rybarczyk \cite{Krzywdzi}, Krishnan
\emph{et al.} \cite{ISIT_RKGRGG}, and we \cite{ZhaoAllerton} present
connectivity results in graph $G(n,K_n,P_n)\bcap G(n,r_n,\mathcal
{A})$. With the network region $\mathcal {A}$ being a square of unit area,
Krzywdzi\'{n}ski and Rybarczyk \cite{Krzywdzi} show that
$G(n,K_n,P_n)\bcap G(n,r_n,\mathcal {A})$ is connected with high
probability if $\pi {r_n}^2 \cdot \frac{{K_n}^2}{P_n} \sim
\frac{c\ln n}{n}$ for any constant $c>8$. Krishnan \emph{et al.}
\cite{ISIT_RKGRGG} improves the condition on $c$ to $c > 2\pi$.
Later we \cite{ZhaoAllerton} derive the critical value $c^*$ of $c$
as
 $\max\{1+\lim_{n \to \infty} \big(\ln \frac{P_n}{{K_n}^2}\big/{\ln n}\big)
,~4\lim_{n \to \infty} \big(\ln \frac{P_n}{{K_n}^2}\big/{\ln n}\big)\}$; namely, %Our
%results represent significantly improved
% conditions for asymptotic connectivity over those of
 graph $G(n,K_n,P_n)\bcap G(n,r_n,\mathcal {A})$
is (resp., is not) connected with high probability for any constant
$c>c^*$ (resp., $c<c^*$). There has not been any analogous result for $k$-connectivity
reported in the literature.

\section{Conclusion and Future Work} \vspace{-2pt}
\label{sec:Conclusion}

In this paper, we consider secure WSNs under the Eschenauer--Gligor (EG)
key predistribution scheme with unreliable links and obtain the
asymptotically exact probability of $k$-connectivity. %Numerical
%simulation is shown to be in agreement with our analytical findings.
A future direction is to consider $k$-connectivity in WSNs employing
 the 
EG scheme under the
disk model \cite{yagan_onoff,ISIT_RKGRGG} in which two nodes have to be within a certain distance
for communication in addition to sharing at least one key.

\iffalse

In this paper, we consider secure WSNs under the EG
key predistribution scheme with unreliable links and rigorously derive the
asymptotically exact probability of $k$-connectivity. %Numerical
%simulation is shown to be in agreement with our analytical findings.
 A future direction is to consider $k$-connectivity in WSNs employing
 the 
EG scheme under the
disk model \cite{yagan_onoff,ISIT_RKGRGG}.
%
%Other key predistribution
%
%Perspectives

\fi

%He
%uses simulation results to conjecture that a similar zero-one law of
%connectivity would be also valid under the disk model, and suggests
%that obtaining such a result analytically would be a significant
%challenge. In this paper, we more than answer this challenge, for
%the following three reasons.
%
%First, we analytically obtain $k$-connectivity (rather than
%1-connectivity) results in key graphs under both the ON/OFF channel
%model \emph{and} the disk model. Second, we consider cases where the
%network area has, or does not have, the boundary effect in the disk
%model.

{ \fontsize{8}{9} \selectfont

\bibliographystyle{abbrv}
\bibliography{related}
  
}

\normalsize

%
%We present several lemmas below which have been used above in the
%paper. We detail the proofs of all the following lemmas in Appendix
%B of our technical report \cite{tech}.

%\section{Useful Lemmas and Their Proofs}
%
%\subsection{Lemmas} \label{subsec:lems}

%
%\begin{lem} \label{lem_psukn}
%
%Consider three distinct nodes $v_i, v_j$ and $v_t$ in graph
%$\mathbb{G}\iffalse_{on}\fi$.
%
%In graph $\mathbb{G}\iffalse_{on}\fi$, for three distinct nodes
%$v_i, v_j$ and $v_t$ in $\mathcal {V}$, let $K_{w v_i }$ (resp.,
%$K_{w v_j }$) be the event that the key rings of nodes $w$ and $v_i$
%(resp., $v_j$) share at least one key. For any given $u = 0, 1, ...,
%K_n$,
%\begin{align}
%\mathbb{P}[(K_{{ wv_i}}\cap K_{ wv_j}) \boldsymbol{\mid} (|S_{i}
%\cap S_{j}| = u)] & \leq
% \frac{ s_n u}{K_n}+ 2{s_n}^2 .
%\nonumber
%\end{align}
%
%
%
%
%\end{lem}

%From (\ref{eqn_epsilonmmmm_m0}),
%\begin{align}
%& \mathbb{P} \big[\mathcal {E} \cap \big( \mathcal {M}_m = \mathcal
%{M}_m^{(0)} \big) \cap \big( \mathcal {T}_m = \mathcal {T}_m^{(0)}
%\big) \big]
%\nonumber  \\
%& \quad =   f\big(n , \mathbb{M}_m^{(0)}\big)
% \mathbb{P}[\mathcal {T}_m = \mathcal {T}_m^{(0)}] \nonumber  \\
%& \quad \quad \times \prod_{i=1}^{m} \mathbb{P}[w \in
%M_{0^{i-1}, 1, 0^{m-i}}^{(0)} \boldsymbol{\mid}%(\mathcal {M}_m = \mathcal
%%{M}_m^{*}) \cap
%\mathcal {T}_m = \mathcal {T}_m^{(0)} ]^{h} \nonumber  \\
%& \quad \quad \times \mathbb{P} [w \in M_{00 \ldots 0}^{(0)}
%\boldsymbol{\mid}%(\mathcal {M}_m = \mathcal
%%{M}_m^{*}) \cap
%\mathcal {T}_m = \mathcal {T}_m^{(0)} ]^{n-hm}. \label{eqn_m00_t00}
%\end{align}

\section*{Appendix} 

 \subsection{Useful Lemmas} 
 
 We present below Lemmas \ref{lem_evalprob} and \ref{lem_psukn}, which are proved in the next subsections. Lemma \ref{lem_evalprob} is used in establishing Propositions \ref{PROP_ONE} and \ref{PROP_SND} in Section \ref{prfseclem1}. The condition $P_n \geq 3K_n $ in Lemma \ref{lem_evalprob} follows for all $n$ sufficiently large given ${K_n}/{P_n} = o(1)$ in Propositions\hspace{3.5pt}\ref{PROP_ONE}\hspace{3.5pt}and\hspace{3.5pt}\ref{PROP_SND}.\hspace{3.5pt}Lemma\hspace{3.5pt}\ref{lem_psukn}\hspace{3.5pt}is\hspace{3.5pt}used\hspace{3.5pt}in\hspace{3.5pt}proving\hspace{3.5pt}Lemma\hspace{3.5pt}\ref{lem_evalprob}.

\begin{lem} \label{lem_evalprob}
 Given $P_n \geq 3K_n $ and any \vspace{1pt}
$\mathcal {T}_m^{*} = (S_1^{*} , S_2^{*}  , \ldots, S_m^{*} )$, with $S_{ij}^{*}$ denoting $S_{i}^{*} \mathlarger{\cap} S_{j}^{*}$, for any node $w \in \overline{\mathcal {V}_m} $,
we \vspace{-2pt} obtain
\begin{align}
 &  \mathbb{P} [w \in M_{0^m} \boldsymbol{\mid} \mathcal {T}_m =
\mathcal {T}_m^{*} ] \geq 1 - m q_n, \textrm{~~~and} \vspace{-2pt} \label{eq_evalprob_3} \\
&\mathbb{P} [w \in M_{0^m} \boldsymbol{\mid} \mathcal
{T}_m = \mathcal {T}_m^{*} ]  \vspace{-2pt} ~~~~ \nonumber \\
& \quad \leq e^{- m q_n + m^2 {q_n}^2 +
  {K_n}^{-1} q_n p_n \sum_{1\leq i <j \leq m}|S_{ij}^{*} |} ; \vspace{-2pt}
   \label{eq_evalprob_4}
\end{align}
and for any $i = 1,2,\ldots,m $, we have \vspace{-2pt}
\begin{align}
 & \mathbb{P}\big[w \in M_{0^{i-1}, 1, 0^{m-i}} \boldsymbol{\mid}
\mathcal {T}_m = \mathcal {T}_m^{*} \big] \leq q_n \vspace{-2pt} \label{eq_evalprob_1} , \textrm{~~~and} \\
&\textstyle{  \prod_{i=1}^{m} \big\{ \mathbb{P}\big[w \in M_{0^{i-1}, 1, 0^{m-i}} 
\boldsymbol{\mid} \mathcal {T}_m = \mathcal {T}_m^{*} \big] \big\}^h} \vspace{-2pt}
\nonumber  \\ 
& \quad  \geq {q_n}^{hm} \textstyle{ \big( 1 - 2 h m^2 {q_n}  - \frac{2 h p_n}{K_n}
\sum _{1\leq i <j \leq m}  |S_{ij}^{*}|\big)}.\vspace{-2pt}\label{hm2qnpnleq}
\end{align}
\iffalse
& \mathbb{P}\big[w \in M_{0^{i-1}, 1, 0^{m-i}}
\boldsymbol{\mid}
\mathcal {T}_m = \mathcal {T}_m^{*} \big]  \nonumber   \\
&  \geq q_n \textstyle{ \big( \hspace{-1pt} 1 \hspace{-2pt} - \hspace{-2pt}
2m{q_n} \hspace{-2pt} - \hspace{-2pt} {K_n}^{-1} p_n \hspace{-2pt}
\sum_{j\in\{1,2,\ldots,m\} \setminus\{i\}} |S_{ij}^{*}|
\hspace{-1pt} \big)}, \label{eq_evalprob_2}

with $\mathcal
{T}_m^{*} = (S_1^{*} , S_2^{*}  , \ldots, S_m^{*} ) $ and
$S_{ij}^{*} = S_{i}^{*} \cap S_{j}^{*}$, we get
\begin{align}
& \prod_{i=1}^{m} \Big\{ \mathbb{P}\big[w \in M_{0^{i-1}, 1, 0^{m-i}} 
\boldsymbol{\mid} \mathcal {T}_m = \mathcal {T}_m^{*} \big] \Big\}^h
\nonumber  \\
& \geq { {q_n}^{hm} \hspace{-1pt} \prod_{i=1}^{m} \bigg[ 1
\hspace{-1pt} - \hspace{-2pt} \bigg( \hspace{-2pt} 2m{q_n}
\hspace{-1pt}
 + \hspace{-1pt} \frac{p_n}{K_n} \hspace{-2pt}
\sum_{j\in\{1,2,\ldots,m\}\setminus\{i\}}
 |S_{ij}^{*}|\bigg)\hspace{-2pt} \bigg] }^h\nonumber  \\
&  \geq {q_n}^{hm} \bigg( 1 - 2 h m^2 {q_n}  - \frac{2 h p_n}{K_n}
\sum _{1\leq i <j \leq m}  |S_{ij}^{*}|\bigg), 
\end{align}

\fi

\end{lem}

%
%We present Lemmas \ref{eqn_Pn_consth}-\ref{lemma_pt2} below to show
%the conclusion that under certain parameter conditions, for any
%non-negative integer $h$, the probability that the connectivity of
%graph $\mathbb{G}\iffalse_{on}\fi$ is $h$ while the minimum node
%degree of graph $\mathbb{G}\iffalse_{on}\fi$ is greater than $h$
%goes to 0 as $n \to \infty$. This conclusion has been used in
%proving property (b) of Theorem \ref{THM1} (viz., Section
%\ref{sec:prf:prob}).

\begin{lem} \label{lem_psukn}

With $\Gamma_{ij} $ denoting the event that an edge exists between distinct nodes
 $v_i$ and $v_j$ in random key graph $G(n,K_n,P_n)$,
if $P_n \geq 3K_n $, then for three distinct nodes \vspace{1pt} $v_i, v_j$ and
$v_t$, we have
$\mathbb{P}[({\Gamma}_{i t} \cap {\Gamma}_{j t} \boldsymbol{\mid}
(|S_{ij}| = u)]  \leq
 {K_n}^{-1} s_n u + 2{s_n}^2 $ for $u = 0, 1,
\ldots, K_n$.\vspace{-2pt}

\end{lem}

 \subsection{The Proof of Lemma \ref{lem_evalprob}}\vspace{-2pt}

 For any node $w \in \overline{\mathcal {V}_m} $,
event $(w \in M_{0^m})$ equals $\overline{\bigcup_{i=1}^{m}
{E_{wv_i}}}$, where $E_{w v_i}$ is the event that there exists an
edge between nodes $w$ and $v_i$ in $\mathbb{G}$. By a union
bound, L.H.S. of (\ref{eq_evalprob_3}) is at least $1 -
\sum_{i=1}^{m} \mathbb{P}[ E_{wv_i} \boldsymbol{\mid}\mathcal {T}_m
= \mathcal {T}_m^{*} ]  = 1 - m q_n$ so that (\ref{eq_evalprob_3}) is proved. And to prove (\ref{eq_evalprob_4}), by the inclusion--exclusion principle, we
get \vspace{-9pt}
\begin{align}
  \mathbb{P} [w  \hspace{-2pt}\in \hspace{-2pt} M_{0^m}
 \hspace{-2pt}\boldsymbol{\mid}  \hspace{-2pt}\mathcal {T}_m  \hspace{-2pt}= \hspace{-2pt} \mathcal {T}_m^{*} ]  
&  \hspace{-1pt}  \leq \hspace{-1pt}  1 \hspace{-1pt} - \hspace{-1pt} \sum_{i=1}^{m} \mathbb{P}[ E_{wv_i}  \hspace{-2pt}
\boldsymbol{\mid} \hspace{-2pt}\mathcal {T}_m  \hspace{-2pt}=  \hspace{-2pt}\mathcal {T}_m^{*} ]  \nonumber
\\& \quad +  \hspace{-1pt}\sum_{1\leq i < j \leq m}  \hspace{-2pt}\mathbb{P}[
E_{wv_{i }} \hspace{-2pt} \bcap \hspace{-2pt} E_{wv_{j}} \hspace{-2pt} \boldsymbol{\mid} \hspace{-2pt}\mathcal {T}_m \hspace{-2pt} = \hspace{-2pt}
\mathcal {T}_m^{*} ] \nonumber.
 \end{align}
Then we use Lemma \ref{lem_psukn}
to further derive
\begin{align}
  & \mathbb{P} [w \in M_{0^m}
\boldsymbol{\mid} \mathcal {T}_m = \mathcal {T}_m^{*} ]  \nonumber
\\&   ~ \leq  1 - m q_n +  {p_n}^2
\sum_{1\leq i <j \leq m}\big({K_n}^{-1} s_n |S_{i j}^{*}| + 2
{s_n}^2
\big) \nonumber%
%\\& ~ =  1 - m q_n  + m(m-1) {q_n}^2  +
% K_n^{-1} q_n p_n \sum_{1\leq i <j \leq m}|S_{i j}^{*}|
%\nonumber
\\&~ \leq e^{- m q_n + m^2 {q_n}^2 +
 {K_n}^{-1} q_n p_n \sum_{1\leq i <j \leq m}|S_{i j}^{*}|}
,\nonumber \vspace{-2pt}
 \end{align}
 where the last step uses $1+x \leq e^{x}$ for any
real $x$.
 
%\end{spacing}
 
% \begin{spacing}{0.95}

 For any node $w \in \overline{\mathcal {V}_m} $,
event $w \in M_{0^{i-1}, 1, 0^{m-i}} $ means that node $w$ has an edge with node $v_i$, but has no edge with any node in $\mathcal {V}_m \setminus \{v_i\} = \{v_j \boldsymbol{\mid} j \in\{1,2,\ldots,m\} \setminus\{i\}\}$. Then (\ref{eq_evalprob_1}) follows since $\mathbb{P}\big[w \in M_{0^{i-1}, 1, 0^{m-i}} 
\boldsymbol{\mid} \mathcal {T}_m = \mathcal {T}_m^{*} \big] $ is at most $ \mathbb{P}[E_{w v_i} \boldsymbol{\mid} \mathcal {T}_m = \mathcal
{T}_m^{*}] = \mathbb{P}[E_{w v_i} ] = q_n$. where the last step uses 
  the independence between event $E_{w v_i} $ and event $ ( \mathcal {T}_m =
\mathcal {T}_m^{*} )$.

We now demonstrate (\ref{hm2qnpnleq}). From the above, we have \vspace{-2pt}
\begin{align}
&  \hspace{-8pt}   \mathbb{P}\big[w \in M_{0^{i-1}, 1, 0^{m-i}} 
\boldsymbol{\mid} \mathcal {T}_m = \mathcal {T}_m^{*} \big] \vspace{-2pt}
\nonumber  \\
&  \hspace{-9pt} = \hspace{-2pt} \textstyle{ \mathbb{P} [ E_{w v_i} \bcap 
\big( \bigcap_{j\in\{1,2,\ldots,m\} \setminus\{i\} }\overline{E_{w v_j} } \big) \boldsymbol{\mid} \mathcal {T}_m =
\mathcal {T}_m^{*}] }  \vspace{-2pt} \nonumber  \\
&  \hspace{-9pt} = \hspace{-2pt} \textstyle{  \mathbb{P} [ E_{w v_i}] \hspace{-2pt} -\hspace{-2pt}  \mathbb{P} [ E_{w v_i}
\hspace{-2pt} \bcap\hspace{-2pt}  \big( \bigcup_{j\in\{1,2,\ldots,m\} \setminus\{i\} }\hspace{-1pt}  E_{w v_j} \hspace{-1pt} \big ) \hspace{-2pt} \boldsymbol{\mid} \hspace{-2pt}  \mathcal {T}_m \hspace{-2pt}  = \hspace{-2pt} 
\mathcal {T}_m^{*}]  },  \vspace{-2pt} \label{tmstarhsp}
\end{align}
where the last step uses $\mathbb{P} [ E_{w v_i} \boldsymbol{\mid} \mathcal {T}_m =
\mathcal {T}_m^{*}] = \mathbb{P} [ E_{w v_i}]  $ since 
 event $E_{w v_i} $ is independent of event $ ( \mathcal {T}_m =
\mathcal {T}_m^{*} )$.

From (\ref{tmstarhsp}) and
$ \mathbb{P} [ E_{w v_i}]  = q_n$, we obtain\vspace{-2pt}
\begin{align}
& {q_n}^{-1}   \mathbb{P}\big[w \in M_{0^{i-1}, 1, 0^{m-i}} 
\boldsymbol{\mid} \mathcal {T}_m = \mathcal {T}_m^{*} \big]   \vspace{-2pt}
\nonumber  \\
&  = \hspace{-2pt}1 \hspace{-1pt} - \hspace{-1pt} {q_n}^{-1}   \textstyle{ \mathbb{P} [   E_{w v_i} \hspace{-1pt}
\bcap \hspace{-1pt} \big(  \bigcup_{j\in\{1,2,\ldots,m\} \setminus\{i\} } E_{w v_j}  \big ) \hspace{-1pt}\boldsymbol{\mid} \hspace{-1pt}\mathcal {T}_m \hspace{-1pt}=\hspace{-1pt}
\mathcal {T}_m^{*}]  }, %\textrm {\hspace{2pt}so\hspace{2pt}that} 
 \vspace{-2pt}\nonumber
\end{align}
%\begin{align}
%& {q_n}^{-1}   \mathbb{P}\big[w \in M_{0^{i-1}, 1, 0^{m-i}} 
%\boldsymbol{\mid} \mathcal {T}_m = \mathcal {T}_m^{*} \big]   \vspace{-2pt}
%\nonumber  \\
%& \hspace{-2pt}= \hspace{-2pt}1 \hspace{-2pt} - \hspace{-2pt} {q_n}^{-1}   \textstyle{ \mathbb{P} [ \hspace{1pt}  E_{w v_i} \hspace{-2pt}
%\bcap \hspace{-2pt} \big( \hspace{-1pt} \bigcup_{j\in\{1,2,\ldots,m\} \setminus\{i\} } E_{w v_j} \hspace{-1pt}\big ) \hspace{-2pt}\boldsymbol{\mid} \hspace{-2pt}\mathcal {T}_m \hspace{-1pt}=\hspace{-1pt}
%\mathcal {T}_m^{*}]  }, %\textrm {\hspace{2pt}so\hspace{2pt}that} 
% \vspace{-2pt}\nonumber
%\end{align}
so that
%which further  induces that
\begin{align}
& \textstyle{  {q_n}^{-hm}  \cdot \prod_{i=1}^{m} \Big\{ \mathbb{P}\big[w \in M_{0^{i-1}, 1, 0^{m-i}}^{(0)}
\boldsymbol{\mid} \mathcal {T}_m = \mathcal {T}_m^{*} \big] \Big\}^h} \vspace{-2pt}
\nonumber  \\
& =  \hspace{-2pt}  \prod_{i=1}^{m}   \textstyle{  \hspace{-2pt}  \Big\{ 1  \hspace{-2pt}  -  \hspace{-2pt} {q_n}^{-1}   \textstyle{ \mathbb{P} [ \hspace{1pt}  E_{w v_i}
 \hspace{-2pt} \bcap \hspace{-2pt}  \big( \bigcup_{j\in\{1,2,\ldots,m\} \setminus\{i\} } E_{w v_j} \big )  \hspace{-2pt}  \boldsymbol{\mid}  \hspace{-2pt}  \mathcal {T}_m  \hspace{-2pt} = \hspace{-2pt} 
\mathcal {T}_m^{*}]  }\Big\}^h}  \vspace{-2pt}\nonumber  \\
& \geq  \hspace{-2pt}  1  \vspace{-2pt} \hspace{-2pt} -   \hspace{-2pt}  h\hspace{-1pt} \sum_{i=1}^{m}   \hspace{-1.5pt} \Big\{  \hspace{-1pt}  {q_n}^{-1}  \textstyle{ \mathbb{P} [  E_{w v_i}
 \hspace{-2.5pt} \bcap  \hspace{-2.5pt}  \big( \bigcup_{j\in\{1,2,\ldots,m\} \setminus\{i\} }  \hspace{-1.5pt}  E_{w v_j} \hspace{-1pt} \big) \hspace{-2.5pt} \boldsymbol{\mid} \hspace{-2.5pt}  \mathcal {T}_m \hspace{-2pt} =\hspace{-2pt} 
\mathcal {T}_m^{*}]  } \hspace{-1pt} \Big\} , \label{ewvjtmstartmq}
\end{align}
where the last step uses the following inequality easily proved by mathematical induction: $\prod_{\ell=1}^{r} (1-x_{\ell}) \geq 1-\sum_{\ell=1}^r x_{\ell} $ for any positive integer $r$ and any $ x_{\ell}$ with $0\leq  x_{\ell}\leq 1$ for $\ell=1,2,\ldots, r$ (we set $r = mh$, with the $mh$ number of $x_l$ as $m$ groups, where the group $i$ for $i=1,2,\ldots, m$ has $m$ members all being $ {q_n}^{-1}    \mathbb{P} [ \hspace{1pt}  E_{w v_i}
 \hspace{-2pt} \bcap \hspace{-2pt}  \big( \bigcup_{j\in\{1,2,\ldots,m\} \setminus\{i\} } E_{w v_j} \big )  \hspace{-2pt}  \boldsymbol{\mid}  \hspace{-2pt}  \mathcal {T}_m  \hspace{-2pt} = \hspace{-2pt} 
\mathcal {T}_m^{*}] $.)

\noindent To  analyze (\ref{ewvjtmstartmq}), we use the union bound and Lemma \ref{lem_psukn} to \vspace{-2pt} get
\begin{align}
&  \textstyle{ \mathbb{P} [ \hspace{1pt}  E_{w v_i}
\bcap \big( \bigcup_{j\in\{1,2,\ldots,m\} \setminus\{i\} } E_{w v_j} \big )\boldsymbol{\mid} \mathcal {T}_m =
\mathcal {T}_m^{*}]  } \vspace{-2pt}
\nonumber  \\
& \quad  \leq  \textstyle{ \sum_{j\in\{1,2,\ldots,m\}\setminus\{i\}}
\mathbb{P} [ E_{w v_i} \cap E_{w v_j} \boldsymbol{\mid} \mathcal
{T}_m = \mathcal {T}_m^{*} ] } \vspace{-2pt} \nonumber  \\
& \quad \leq  \textstyle{  \sum_{j\in\{1,2,\ldots,m\}\setminus\{i\}}{p_n}^2 \big(
{K_n}^{-1} s_n  |S_{ij}^{*}| + 2{s_n}^2 \big) } \vspace{-2pt}  \nonumber  \\
& \quad \leq  \textstyle{ 2m {q_n}^2  + {K_n}^{-1} p_n q_n
\sum_{j\in\{1,2,\ldots,m\}\setminus\{i\}}   |S_{ij}^{*}| }, \vspace{-2pt} \nonumber
\end{align}
which is substituted into (\ref{ewvjtmstartmq}) to establish (\ref{hm2qnpnleq})  by \vspace{-2pt}
\begin{align}
& \textstyle{  {q_n}^{-hm}  \cdot \prod_{i=1}^{m} \big\{ \mathbb{P}\big[w \in M_{0^{i-1}, 1, 0^{m-i}}
\boldsymbol{\mid} \mathcal {T}_m = \mathcal {T}_m^{*} \big] \big\}^h} \vspace{-2pt}
\nonumber  \\
&   \geq   1 - \textstyle{ h \sum_{i=1}^{m}  \vspace{-2pt} \big\{ 2m {q_n}  +  \textstyle{ {K_n}^{-1} p_n
\sum_{j\in\{1,2,\ldots,m\}\setminus\{i\}}   |S_{ij}^{*}| } \big\} } \nonumber  \\
&   \geq  1 - 2 h m^2 {q_n}  - \textstyle{ \frac{2 h p_n}{K_n}
\sum _{1\leq i <j \leq m}  |S_{ij}^{*}|}. \vspace{-2pt}
\end{align}

 \subsection{The Proof of Lemma \ref{lem_psukn}} \vspace{-2pt}

We use the inclusion--exclusion principle to obtain\vspace{-2pt}
\begin{align}
& \mathbb{P}[\Gamma_{{ it}} \hspace{1.5pt}\mathlarger{\cap}  \hspace{1.5pt}\Gamma_{ jt} \boldsymbol{\mid}
(|S_{ij}| = u)]  \vspace{-2pt} \nonumber \\   & ~= \mathbb{P}[\Gamma_{{it}}
\boldsymbol{\mid} (|S_{ij}| = u)] + \mathbb{P}[
\Gamma_{jt} \boldsymbol{\mid} (|S_{ij}| = u)]  \vspace{-2pt}\nonumber \\
& ~\quad- \mathbb{P}[ {\Gamma_{{ it}}} \hspace{1.5pt}\mathlarger{\cup} \hspace{1.5pt}
 {\Gamma_{ jt}}\boldsymbol{\mid} (|S_{ij}| = u)]
\vspace{-2pt}
 \nonumber \\  &~ = \textstyle{ 2 s_n - 1 + \binom{P_n -
(2K_n -u)}{K_n}\big/\binom{P_n}{K_n} }\vspace{-2pt}, \label{sn1Pnsiju2Knuan}
\end{align}
in view that event $(|S_{ij}| = u)$ is independent of each of $\Gamma_{{ it}} $ and $ \Gamma_{ jt} $, and event ${\Gamma_{{ it}}} \hspace{1.5pt}\mathlarger{\cup} \hspace{1.5pt}
 {\Gamma_{ jt}} $ means  $S_t  \hspace{1.5pt}\mathlarger{\cap} \hspace{1.5pt} (S_i  \hspace{1.5pt}\mathlarger{\cup}  \hspace{1.5pt}S_j) \neq \emptyset$. 

By \cite[Lemma 5.1]{yagan_onoff} and \cite[Fact 2]{ZhaoYaganGligor}, we derive
\begin{align}
(1-s_n)^{\frac{2K_n -u}{K_n}}  & \leq \vspace{-2pt}1 - \textstyle{\frac{s_n(2K_n -u)}{K_n} +\frac{1}{2}\big(\frac{s_n(2K_n -u)}{K_n}\big)^2}  \nonumber \\& \leq 1- 2s_n + {K_n}^{-1} s_n u + 2{s_n}^2 ,\vspace{-2pt} \nonumber
\end{align}
which is substituted into (\ref{sn1Pnsiju2Knuan}) to complete the proof.  

\end{spacing}

\end{document}